\documentclass[11pt,reqno,twoside]{amsart}
\usepackage{amsmath,amssymb,amsfonts,amsthm,mathtools}
\usepackage{geometry}
\usepackage{indentfirst}
\usepackage{enumitem}
\usepackage{needspace}
\usepackage[noadjust]{cite}
\usepackage[colorlinks=true,citecolor=red,linkcolor=blue]{hyperref}
\allowdisplaybreaks[4]
\makeatletter
\@namedef{subjclassname@2020}{\textup{2020} Mathematics Subject Classification}
\makeatother
\hypersetup{
 pdftitle={Localization for Lipschitz monotone quasi-periodic operators with power-law hopping},
 pdfauthor={Hongyi Cao; Yunfeng Shi; Peizhen Wang},
 pdfsubject={Localization; Lipschitz monotone potentials; quasi-periodic operators; power-law hopping}
}
\numberwithin{equation}{section}
\theoremstyle{plain}
\newtheorem{theorem}{Theorem}[section]
\newtheorem{proposition}[theorem]{Proposition}
\newtheorem{lemma}[theorem]{Lemma}
\newtheorem{corollary}[theorem]{Corollary}
\theoremstyle{definition}

\newtheorem{remark}[theorem]{Remark}

\newcommand{\dist}{\operatorname{dist}}
\newcommand{\Spec}{\operatorname{Spec}}
\newcommand{\eps}{\varepsilon}
\newcommand{\om}{\omega}
\newcommand{\Sch}{\operatorname{Sch}}
\newcommand{\Res}{\operatorname{Res}}
\newcommand{\Id}{\mathcal I}
\newcommand{\fZ}{Z}

\newcommand{\cW}{\mathcal W}

\newcommand{\fB}{B}
\newcommand{\fF}{F}
\newcommand{\fJ}{J}
\newcommand{\fL}{\Lambda}
\newcommand{\fO}{\Omega}
\newcommand{\fP}{P}
\newcommand{\fQ}{Q}
\newcommand{\fS}{S}
\newcommand{\fT}{\Theta}
\newcommand{\fX}{X}
\newcommand{\fY}{Y}
\newcommand{\dc}{\mathrm{DC}}
\newcommand{\matnorm}[2]{\lVert #1\rVert_{#2}}

\newcommand{\cA}{\mathcal A}
\newcommand{\cM}{\mathcal M}
\newcommand{\cD}{\mathcal D}
\newcommand{\cH}{\mathcal H}
\newcommand{\cT}{\mathcal T}
\newcommand{\cU}{\mathcal U}
\title[Localization for monotone operators with power-law hopping]{Localization for Lipschitz monotone quasi-periodic operators with power-law hopping}

\author[Cao]{Hongyi Cao}
\address[Cao]{School of Mathematical Sciences, Peking University, Beijing 100871, China}
\email{chyyy@stu.pku.edu.cn}

\author[Shi]{Yunfeng Shi}
\address[Shi]{School of Mathematics, Sichuan University, Chengdu 610064, China}
\email{yunfengshi@scu.edu.cn}

\author[Wang]{Peizhen Wang}
\address[Wang]{School of Mathematics, Sichuan University, Chengdu 610064, China}
\email{pzw@scu.edu.cn}

\date{November 12, 2026}

\begin{document}

\begin{abstract}
We prove spectral and dynamical localization for quasi-periodic
operators on $\mathbb Z^d$ with Lipschitz monotone potentials and
power-law hopping.
For sufficiently small coupling, localization holds at every
allowed phase whenever the hopping decay exponent $r$ exceeds
the Diophantine exponent $\tau$.
The proof is based on multiscale analysis with exponentially growing
scales. We use homotopy in the coupling constant to select the local
eigenvalue branches and Schur complement estimates to control the
successive elimination of coordinates.
\end{abstract}

\maketitle

\section{Introduction}
\subsection{The model}\label{sec:model}
We study the following family of quasi-periodic operators on
$\ell^2(\mathbb Z^d)$ $(d\ge 1)$:
\begingroup\postdisplaypenalty=100
\begin{equation}\label{2026092501}
 (\cH(\theta) u)(n)
 =f(\theta+n\cdot\om)u(n)
 +\eps\sum_{m\in\mathbb Z^d}\phi(n-m)u(m),
 \qquad n\in\mathbb Z^d.
\end{equation}
\endgroup
Here $\theta\in\mathbb R$ is the phase,
$\eps\in\mathbb R$ is the coupling constant with $|\eps|\ll 1$ and
$\om\in\mathbb R^d$ is the frequency.
The notation ``$\cdot$'' in \eqref{2026092501} denotes the Euclidean inner product
on $\mathbb R^d$.
We assume $\omega\in\dc$, where
\begin{equation}\label{2026092502}
 \dc
 =\left\{\omega\in\mathbb R^d:
   \|k\cdot\omega\|_{\mathbb T}\ge\gamma\langle k\rangle^{-\tau}
   \text{ for all }k\in\mathbb Z^d\setminus\{o\}\right\}
\end{equation}
with fixed $\gamma>0$ and $\tau>d$.

The symbols $o$, $\langle k\rangle$ and $\|\cdot\|_{\mathbb T}$
in \eqref{2026092502} are defined as follows.
For $x=(x_1,\ldots,x_d)\in\mathbb R^d$, set
\begin{equation}\label{2026092503}
 o=(0,\ldots,0)\in\mathbb R^d,\qquad
 |x|=\max_{1\le i\le d}|x_i|,\qquad
 \langle x\rangle=1+|x|.
\end{equation}
For $\alpha\in\mathbb R$, set
\begin{equation}\label{2026092504}
 \mathbb T=\mathbb R/\mathbb Z,\qquad
 \|\alpha\|_{\mathbb T}=\min_{m\in\mathbb Z}|\alpha-m|.
\end{equation}

For some $r>d$, the hopping $\phi:\mathbb Z^d\to\mathbb C$ in \eqref{2026092501} satisfies
\begingroup\postdisplaypenalty=100
\begin{equation}\label{2026092601}
 \phi({o})=0,\qquad \overline{\phi(k)}=\phi(-k),\qquad
 |\phi(k)|\le\mathrm C_0\langle k\rangle^{-r}
 \quad(k\in\mathbb Z^d),
\end{equation}
\endgroup
where
\begin{equation}\label{2026092602}
	\mathrm C_0=\mathrm C_0(r,\phi)
	=\max\left\{1,\sup_{k\in\mathbb Z^d}
	\langle k\rangle^r|\phi(k)|\right\}<\infty.
\end{equation}
For the potential, we assume that $f$ is $1$-periodic,
\begin{equation}\label{2026092603}
 f:\mathbb R\to[-\infty,+\infty),\qquad f(\theta+1)=f(\theta),\qquad
 f(\theta)\in\mathbb R\quad(0<\theta<1),
\end{equation}
and Lipschitz monotone on $[0,1)$: for some $L>0$,
\begin{equation}\label{2026092604}
 f(\theta_2)-f(\theta_1)\ge L(\theta_2-\theta_1),\qquad 0\le \theta_1<\theta_2<1.
\end{equation}
We call $L$ the monotonicity constant.

We consider $\cH(\theta)$ only at phases for which
$f(\theta+n\cdot\om)$ is finite for every $n\in\mathbb Z^d$. Denote
\begin{equation}\label{2026092605}
 \fT=\{\theta\in\mathbb R:
     f(\theta+n\cdot\om)\in\mathbb R
     \text{ for every }n\in\mathbb Z^d\}.
\end{equation}
Thus $\fT=\mathbb R$ if $f(0)\in\mathbb R$, whereas
$\fT=\mathbb R\setminus\{m+n\cdot\om:m\in\mathbb Z,\ n\in\mathbb Z^d\}$
if $f(0)=-\infty$ (equivalently, $f$ is unbounded below).
For every $\theta\in\fT$, $\cH(\theta)$ is self-adjoint.

\subsection{Motivation}\label{sec:motivation}
At $\eps=0$, the operator \eqref{2026092501} has an eigenbasis
supported at single lattice points.
We ask whether a complete system of localized eigenvectors
persists for sufficiently small $\eps$.
For nearest-neighbor hopping, the Maryland model is an early
example (\cite{GrempelFishmanPrange1982}).
In 1983, Bellissard, Lima and Scoppola (\cite{BLS83}) established
perturbative localization for meromorphic monotone potentials
in arbitrary dimension.
The restriction to meromorphic potentials led to the question
of localization under a quantitative monotonicity condition alone.

In 2019, Jitomirskaya and Kachkovskiy
(\cite{JitomirskayaKachkovskiy2019}) introduced a general class
of bounded Lipschitz monotone potentials.
The all-coupling localization theorem applies in one dimension
to Diophantine frequencies and almost every phase.
Higher-dimensional perturbative results were subsequently obtained
for unbounded monotone potentials in 2022
(\cite{KachkovskiyParnovskiShterenberg2022}) and for H\"older monotone
potentials in 2025 (\cite{KPS2025}).
Cao, Shi and Zhang (\cite{CSZ2025}) also established localization
uniformly in phase for bounded and unbounded Lipschitz monotone potentials.

For the hopping term, perturbative localization also holds
under exponential decay (\cite[Remark~1.5]{CSZ2025}).
More recently, Cao (\cite{Cao2026}) obtained arithmetic localization
for one-dimensional exponentially decaying hopping under an
explicit frequency condition and a small-coupling assumption.
For power-law hopping, Shi (\cite{Shi2023}) and Shi and Wen
(\cite{ShiWen2023}) proved localization results in 2023.
For the Maryland potential, a sufficient pointwise hopping
decay condition in \cite{Shi2023,ShiWen2023} is $r>\tau+d$.
We ask whether perturbative localization holds for the original
operator with a general Lipschitz monotone potential and
power-law hopping under $r>\tau$.
\subsection{Main results}\label{sec:main-results}
Our main result is as follows.

\begin{theorem}\label{thm:main}
Fix an integer $d\ge1$, numbers $\tau>d$, $\gamma,L>0$ and $r>\tau$.
Assume \eqref{2026092603}, \eqref{2026092604} and \eqref{2026092601}.
There is a constant
$\mathrm C_{1}=\mathrm C_{1}(d,r,\tau,\gamma,L,\phi)>0$ such that,
for every $0<\rho<r-\tau$, there are constants
$\eps_0=\eps_0(d,\tau,\gamma,L,r,\rho,\phi)>0$ and
$\mathrm C_{2}=\mathrm C_{2}(d,r,\tau,\gamma,L,\rho,\phi)>0$
with the following property.
For every $\om\in\dc$ and $|\eps|<\eps_0$,
there is a real $1$-periodic function $E_\infty:\fT\to\mathbb R$
such that
\begingroup\postdisplaypenalty=10000
\begin{equation*}
 \|E_\infty-f\|_\infty\le\mathrm C_{1}\eps^2,\qquad
 E_\infty(\theta_2)-E_\infty(\theta_1)\ge L(\theta_2-\theta_1)
 \quad(\theta_1<\theta_2,\ \theta_1,\theta_2\in\fT\cap[0,1)),
\end{equation*}
where $\|\cdot\|_\infty$ denotes the supremum norm on
$\fT\cap[0,1)$.
\endgroup
For every $\theta\in\fT$, the operator $\cH(\theta)$ has a
complete orthonormal eigenbasis
$(\psi_{\theta,n})_{n\in\mathbb Z^d}$ satisfying
\[
 \cH(\theta)\psi_{\theta,n}
 =E_\infty(\theta+n\cdot\om)\psi_{\theta,n}.
\]
For $m\ne n$,
$E_\infty(\theta+m\cdot\om)\ne E_\infty(\theta+n\cdot\om)$.
The vectors $\psi_{\theta,n}$ can be chosen so that
\begingroup\postdisplaypenalty=10000
\begin{equation}\label{2026092606}
 \sup_{n\in\mathbb Z^d}
 \sum_{m\in\mathbb Z^d}
 |\psi_{\theta,n}(m)-\delta_{mn}|
       \langle m-n\rangle^\rho
 \le\mathrm C_{2}|\eps|,
\end{equation}
where $\delta_{mn}=1$ if $m=n$ and $\delta_{mn}=0$ otherwise.
\endgroup
Moreover, for the unitary operator
$\exp(-\mathrm i u\cH(\theta))$ defined by the spectral theorem,
\label{def:unitary-evolution}
\begingroup\postdisplaypenalty=10000
\begin{equation}\label{2026092607}
 \sup_{u\in\mathbb R}\sup_{n\in\mathbb Z^d}
 \sum_{m\in\mathbb Z^d}
 |\langle e_m,\exp\left(-\mathrm i u\cH(\theta)\right)e_n\rangle|
       \langle m-n\rangle^\rho
 \le 3,
\end{equation}
where $\langle\cdot,\cdot\rangle$ is the $\ell^2$ inner product,
linear in the second argument, and
$e_n(m)=\delta_{nm}$.
\endgroup
\end{theorem}

\begin{remark}\label{rem:potential-class}
Theorem~\ref{thm:main} applies to power-law hopping and to
both bounded and unbounded potentials satisfying
\eqref{2026092603} and \eqref{2026092604}.
No meromorphic extension of $f$ is assumed, whereas
meromorphic continuation is used in \cite{BLS83,ShiWen2023}.
\end{remark}

\begin{remark}\label{rem:hopping-exponent}
The sufficient hopping decay condition in Theorem~\ref{thm:main}
is $r>\tau$.
For the Maryland potential, the sufficient pointwise condition is
$r>\tau+d$ in \cite[Remark~2.3 and Corollary~2.10]{Shi2023}
and \cite[Remark~1.5]{ShiWen2023}.
On this common class of operators, Theorem~\ref{thm:main}
removes the additional $d$ from the sufficient decay condition.
Theorem~\ref{thm:main} also applies to the general Lipschitz monotone
class in \eqref{2026092603}--\eqref{2026092604}.
The endpoint $r=\tau$ is not covered.
\end{remark}

\begin{remark}\label{rem:uniform-localization}
Theorem~\ref{thm:main} gives pure point spectrum at every
$\theta\in\fT$, with no additional arithmetic restriction on the phase.
The threshold $\eps_0$ and the bounds
\eqref{2026092606}--\eqref{2026092607} are independent of
$\theta\in\fT$.
For every fixed $0<\rho<r-\tau$, \eqref{2026092606} implies
\[
 |\psi_{\theta,n}(m)|\le\mathrm C_{2}|\eps|
       \langle m-n\rangle^{-\rho},\qquad m\ne n.
\]
Moreover, \eqref{2026092607} gives polynomial dynamical
localization with exponent $\rho$.
\end{remark}

\subsection{Related work}\label{sec:related-work}
The results in Subsection~\ref{sec:motivation} arise from two
different approaches to localization.
Early perturbative constructions for almost-periodic operators
were developed by Craig and P\"oschel in 1983
(\cite{Craig1983,Poschel1983}).
For monotone potentials, the KAM proof of Bellissard, Lima and
Scoppola (\cite{BLS83}) uses meromorphic continuation to control
the small denominators.
The one-dimensional nonperturbative approach instead relates
localization to the Lyapunov exponent.
Jitomirskaya and Kachkovskiy (\cite{JitomirskayaKachkovskiy2019})
combined monotonicity estimates with large deviation estimates
for finite-volume determinants.
This provided a way to treat potentials without the special
structure of the Maryland potential.

The nonperturbative theory was extended in several directions.
In 2019, Kachkovskiy (\cite{Kachkovskiy2019}) treated unbounded
monotone potentials under a logarithmic integrability condition,
which ensures that the Lyapunov exponent is finite.
The upper Lipschitz bound in the bounded case was also removed
(\cite[Remark~1.4]{Kachkovskiy2019}).
Meanwhile, arithmetic localization studies the competition
between the Lyapunov exponent and rational approximation of
the frequency.
Jitomirskaya and Liu (\cite{JitomirskayaLiu2017}) determined the
arithmetic spectral transitions for the Maryland model in 2017.
Jitomirskaya and Kachkovskiy (\cite{JitomirskayaKachkovskiy2024})
established the sharp arithmetic localization criterion for
Lipschitz monotone potentials in 2024.
More recently, Jitomirskaya, Levi and Stoiber
(\cite{JitomirskayaLeviStoiber2026}) treated a class of monotone potentials
with flat regions.
The proof in \cite{JitomirskayaLeviStoiber2026} recovers quantitative
variation of finite-volume eigenvalues by estimating eigenfunction
mass outside the flat regions.
For exponentially decaying long-range hopping, Cao
(\cite{Cao2026}) extended the nonperturbative arithmetic method
by estimating Green's functions without a Lyapunov exponent.

In higher dimensions, the recent perturbative proofs use
different constructions of the eigenvalues and eigenvectors.
Kachkovskiy, Parnovski and Shterenberg
(\cite{KachkovskiyParnovskiShterenberg2022}) proved convergence
of perturbation series by grouping terms to exploit cancellations.
In \cite{KPS2025}, Kachkovskiy, Parnovski and Shterenberg used local partial
diagonalization and finite-volume eigenvalue ordering.
The ordered branches recover the monotonicity needed to
control small denominators throughout the iteration.
Cao, Shi and Zhang (\cite{CSZ2025}) instead constructed
Rellich functions on increasing finite regions and used
Green's function estimates to pass between scales.
Monotonicity separates points whose local eigenvalues are
close to the same energy.
For bounded potentials, a Schur complement comparison across
discontinuities ensures
that the selected Rellich functions preserve Lipschitz monotonicity.

For power-law hopping, the estimates must also control sums of
slowly decaying matrix entries.
Shi (\cite{Shi2023}) developed a Nash--Moser reducibility scheme
in Sobolev matrix norms.
Tame product estimates and smoothing compensate for the loss of
regularity in the homological equations.
The general theorem gives reducibility after a small correction
of the diagonal potential
(\cite[Theorem~2.1]{Shi2023}).
The theorem for the original operator assumes that the small-denominator
bounds persist under all sufficiently small diagonal perturbations
in the chosen Banach algebra
(\cite[Theorem~2.2]{Shi2023}).
Shi and Wen (\cite{ShiWen2023}) used meromorphic monotonicity
to obtain diagonalization uniformly in phase.
In the present paper, we combine the Rellich-function construction
with Schur complement estimates adapted to power-law decay.

\subsection{Main ideas and new ingredients of the proof}\label{sec:main-ideas}
We use the Rellich-function construction of Cao, Shi and Zhang
(\cite{CSZ2025}) as the starting point.
The main changes are as follows.

\smallskip\noindent
\textbf{1. Monotonicity of the local eigenvalues.}
As in \cite{CSZ2025}, we construct local eigenvalues on
increasing finite regions and use monotonicity to separate
the resonant points.
The finite-volume ordering argument in
\cite[Proposition~4.1 and Remark~4.2]{KPS2025} gives eigenvalue
branches with the same Lipschitz monotonicity as the potential.
By continuation from zero coupling to $\eps$, we identify each
isolated local eigenvalue with the ordered branch associated
with the diagonal entry $f(\theta)$
(Proposition~\ref{prop:homotopy-monotonicity}).
The branch identification gives the required monotonicity in place
of the comparison across jumps used for bounded potentials
in \cite{CSZ2025}.

\smallskip\noindent
\textbf{2. Schur estimates under the condition $r>\tau$.}
Shi's Nash--Moser iteration controls the conjugating matrices
in Sobolev norms (\cite{Shi2023}).
We estimate the entries of each Schur complement directly and
use the separation of the sites retained from earlier steps
to bound the sums over coordinates being eliminated.
Once the diagonal entries are bounded away from zero,
Lemma~\ref{lem:sparse-schur} gives inverse estimates and preserves
the off-diagonal decay exponent $r$.
After elimination, however, the diagonal entries depend on the finite region.
Corollary~\ref{cor:diagonal-comparison} bounds the change of a
diagonal entry when the region is enlarged.
We can therefore use local eigenvalues to obtain the diagonal
lower bounds needed for elimination in the larger region.
The estimates close under $r>\tau$ for the prescribed Lipschitz
monotone potential.

\smallskip\noindent
\textbf{3. Exponentially growing multiscale lengths.}
In \cite[Definition~2.15]{CSZ2025}, the spatial scale is squared
at each induction step.
In our construction, the spatial scale increases by the fixed
factor $16$.
Equation~\eqref{2026100701} controls the accumulation of the Schur
estimates by a geometric series with ratio $16^{\tau-r}<1$.
The summability and the uniform region comparison
\eqref{2026101305} close the induction with a coupling
threshold independent of the step.

\subsection{Outline of this paper}\label{sec:proof-outline}
Section~\ref{sec:preliminaries} collects notation and standard matrix facts.
Section~\ref{sec:initial-construction} constructs the first local eigenvalue
and the corresponding eigenvector, then continues with the remaining equations.
Section~\ref{sec:inductive-hypotheses} states the induction hypotheses.
Section~\ref{sec:induction-step} proves the region comparison, performs
the next elimination and estimates the changes of the local eigenvalues
and corresponding eigenvectors for every $n\ge1$.
Section~\ref{sec:localization} passes to the limit and proves completeness
and dynamical localization.
Appendix~\ref{app:matrix-estimates} contains the sparse matrix estimates;
Appendix~\ref{app:technical-verifications} gives the geometric and
domain verifications.

\Needspace{7\baselineskip}
\section{Preliminaries}\label{sec:preliminaries}

\subsection{Notation}\label{sec:notation}
\begin{enumerate}[label=\arabic*.,ref=\arabic*,leftmargin=*,itemsep=\medskipamount]
\item \textbf{Distances and phases.}
With $|\cdot|$ as in \eqref{2026092503}, define
\begin{equation}\label{2026092701}
 \dist(x,S)=\inf_{y\in S}|x-y|,
 \qquad x\in\mathbb R^d,\quad S\subset\mathbb R^d,
\end{equation}
with $\dist(x,\varnothing)=+\infty$.
For $\alpha\in\mathbb R$, denote the fractional part by
\begin{equation}\label{2026092702}
 \{\alpha\}=\alpha-\lfloor\alpha\rfloor\in[0,1).
\end{equation}
Unless another domain is specified, phase-dependent functions are
considered on the set $\fT$ in \eqref{2026092605}.

\item \textbf{Lattice sets.}
For $\fZ,\fB\subset\mathbb Z^d$ and $p\in\mathbb Z^d$, denote
\begin{equation}\label{2026092703}
 \fZ^c=\mathbb Z^d\setminus\fZ,\qquad
 p+\fZ=\{p+x:x\in\fZ\},\qquad
 \fZ+\fB=\{x+y:x\in\fZ,\ y\in\fB\}.
\end{equation}
The cardinality of a finite set $\fZ$ is denoted by $\#\fZ$.
If $o\in\fZ$, write
$\fZ^\circ=\fZ\setminus\{o\}$.\label{def:punctured-set}
For a set $B\subset\mathbb Z^d$ containing $o$ and $p\in\mathbb Z^d$,
denote
\begin{equation}\label{2026092704}
 B(p)=p+B,\qquad B^\circ(p)=B(p)\setminus\{p\}.
\end{equation}
For $R\ge0$, let
\begin{equation}\label{2026092705}
 \fL_R=\{n\in\mathbb Z^d:|n|\le R\}.
\end{equation}
For $M>0$, a set $\fX\subset\mathbb Z^d$ is
$M$-separated if $|p-q|\ge M$ for all distinct $p,q\in\fX$.
\label{def:separated-set}

\item \textbf{Vectors and matrices.}\label{item:matrix-notation}
For a vector $v$ whose index set contains $\fZ$, write
$\left.v\right|_{\fZ}=(v(n))_{n\in\fZ}$. For $1\le p<\infty$, let
$\|v\|_{\ell^p(\fZ)}=(\sum_{n\in\fZ}|v(n)|^p)^{1/p}$;
we write $\|v\|_p$ when the index set is fixed.
For an operator $\cM$, write $\cM(n,m)=\langle e_n,\cM e_m\rangle$.
The adjoint is denoted by $\cM^*$ and $\|\cM\|_{\rm op}$ denotes
the $\ell^2$ operator norm. For a finite matrix,
$\cM^*(n,m)=\overline{\cM(m,n)}$; for a column vector $v$,
the symbol $v^*$ denotes the conjugate transpose of $v$.
Let $\Id_\fZ$ be the identity on $\ell^2(\fZ)$ and let
$\Id=\Id_{\mathbb Z^d}$.
For $p\in\mathbb Z^d$, define the translation by
\begin{equation}\label{2026092706}
 (\cU(p)v)(n)=v(n+p).
\end{equation}
For a matrix $\cM$, let
$\cM_{\fZ,\fB}=(\cM(p,q))_{p\in\fZ,\,q\in\fB}$ and write
$\cM_\fZ=\cM_{\fZ,\fZ}$.
A point $p$ in a submatrix index denotes the singleton $\{p\}$; thus
$\cM_{p,\fB}=\cM_{\{p\},\fB}$ and
$\cM_{\fZ,q}=\cM_{\fZ,\{q\}}$.
Phase and energy variables are written as arguments.

\end{enumerate}

\subsection{Operators and finite restrictions}\label{sec:operator-preliminaries}
For $q=1,2$ and a real family $(v_n)_{n\in Z}$, the symbol
$\operatorname{diag}_{n\in Z}v_n$ denotes the multiplication operator
on $\ell^q(Z)$ with the natural domain.
The multiplication operator is closed, and finitely supported
vectors form a core. On $\ell^2(Z)$ the operator is self-adjoint.

The hopping operator is defined by
\begin{equation}\label{2026092707}
 (\cT v)(n)=\sum_{m\in\mathbb Z^d}\phi(n-m)v(m).
\end{equation}
By \eqref{2026092601} and $r>d$, $\phi\in\ell^1(\mathbb Z^d)$;
hence $\cT$ is bounded on $\ell^1(\mathbb Z^d)$ and
$\ell^2(\mathbb Z^d)$.
The symmetry in \eqref{2026092601} implies that $\cT$ is
self-adjoint on $\ell^2(\mathbb Z^d)$.

The phase set $\fT$ satisfies
$\fT+q=\fT$ and $\fT+p\cdot\om=\fT$ for $q\in\mathbb Z$ and
$p\in\mathbb Z^d$. Direct substitution into \eqref{2026092501}
proves
\begin{equation}
 \cH(\theta+1)=\cH(\theta),\qquad \cH(\theta+p\cdot\om)\cU(p)=\cU(p)\cH(\theta),
 \qquad \theta\in\fT,\quad p\in\mathbb Z^d.
 \label{2026092708}
\end{equation}
Here $\cU(p)$ is defined in \eqref{2026092706}.

For $\fZ\subset\mathbb Z^d$, let
$\Pi_\fZ:\ell^2(\mathbb Z^d)\to\ell^2(\fZ)$ be restriction to
$\fZ$. The adjoint $\Pi_\fZ^*$ extends a vector by zero outside
$\fZ$. For finite $\fZ$, $\theta\in\fT$ and $z\in\mathbb C$, let
\begin{align}
 \cH_\fZ(\theta)&=\Pi_\fZ\cH(\theta)\Pi_\fZ^*,\qquad
 \cT_\fZ=\Pi_\fZ\cT\Pi_\fZ^*,\label{2026092801}\\
 \cA_\fZ(z;\theta)&=z\Id_\fZ-\cH_\fZ(\theta).
 \label{2026092802}
\end{align}
When $\cA_\fZ(z;\theta)$ is invertible,
$\cA_\fZ(z;\theta)^{-1}$ is the finite-volume Green's function.

\subsection{Schur complements}\label{sec:schur-preliminaries}
Let $\cM$ be a matrix on a finite disjoint union $\fP\sqcup\fQ$
and assume that $\cM_\fQ$ is invertible.
The Schur complement after eliminating $\fQ$ is
\begin{equation}\label{2026092803}
 \Sch_\fQ(\cM)
 =\cM_\fP-\cM_{\fP,\fQ}\cM_\fQ^{-1}\cM_{\fQ,\fP}.
\end{equation}
For the decomposition $\fP\sqcup\fQ$, define
$\cW(\cM):\mathbb C^\fP\to\mathbb C^{\fP\sqcup\fQ}$ by
\begin{equation}\label{2026092804}
 \cW(\cM)=
 \begin{pmatrix}\Id_\fP\\-\cM_\fQ^{-1}\cM_{\fQ,\fP}\end{pmatrix},
\end{equation}
with row blocks in the order $\fP,\fQ$ and column set $\fP$.
Solving the equations with indices in $\fQ$, we obtain
\begin{align*}
 \cM v=0
 &\quad\Longleftrightarrow\quad
 \left.v\right|_{\fQ}=-\cM_\fQ^{-1}\cM_{\fQ,\fP}\left.v\right|_{\fP},\qquad
 \Sch_\fQ(\cM)\left.v\right|_{\fP}=0,\\
 \ker\cM&=\cW(\cM)\ker\Sch_\fQ(\cM).
\end{align*}
Block multiplication also proves
\begin{equation*}
 \cM\cW(\cM)=
 \begin{pmatrix}\Sch_\fQ(\cM)\\0\end{pmatrix},\qquad \det\cM=\det\cM_\fQ\det\Sch_\fQ(\cM).
\end{equation*}
If $\cM=\cM^*$, then
$\cW(\cM)^*\cM\cW(\cM)=\Sch_\fQ(\cM)$.

Let $\cM$ act on a finite disjoint union $\fP\sqcup\fQ\sqcup U$.
If $\cM_\fQ$ and $(\Sch_\fQ(\cM))_U$ are invertible, then
$\cM_{\fQ\cup U}$ is invertible and
\begin{equation}\label{2026092805}
 \Sch_U(\Sch_\fQ(\cM))=\Sch_{\fQ\cup U}(\cM).
\end{equation}
Indeed, block multiplication gives
$\det\cM_{\fQ\cup U}
=\det\cM_\fQ\det(\Sch_\fQ(\cM))_U\ne0$,
so $\cM_{\fQ\cup U}$ is invertible.
For each prescribed $\left.v\right|_{\fP}$, solving first for $\left.v\right|_{\fQ}$ and then
for $\left.v\right|_{U}$ produces the unique solution of the equations with
indices in $\fQ\cup U$. Substitution into the equations indexed by
$\fP$ proves \eqref{2026092805}.
For $\fQ=\varnothing$, sums over $\fQ$ are zero,
$\Id_\varnothing$ and $\cM_\varnothing^{-1}$ denote the unique operator
on $\mathbb C^\varnothing$, and $\Sch_\varnothing(\cM)=\cM$.

For Hermitian matrices $\cM,\widetilde{\cM}$ on the same
finite set $\fP\sqcup\fQ$, assume that $\cM_\fQ$ and
$\widetilde{\cM}_\fQ$ are invertible. Then
\begin{equation}\label{2026092806}
 \Sch_\fQ(\widetilde{\cM})-\Sch_\fQ(\cM)
 =\cW(\cM)^*(\widetilde{\cM}-\cM)\cW(\widetilde{\cM}).
\end{equation}
Here $\cW(\cM)$ and $\cW(\widetilde{\cM})$ are defined by
\eqref{2026092804} with row blocks $\fP,\fQ$.
By \eqref{2026092804} and the Hermitian symmetry of
$\cM$ and $\widetilde{\cM}$, block multiplication gives
\begin{equation*}
 \cW(\cM)^*\cM\cW(\widetilde{\cM})=\Sch_\fQ(\cM),\qquad \cW(\cM)^*\widetilde{\cM}\cW(\widetilde{\cM})
 =\Sch_\fQ(\widetilde{\cM}).
\end{equation*}
Subtracting the identities for
$\cW(\cM)^*\widetilde{\cM}\cW(\widetilde{\cM})$ and
$\cW(\cM)^*\cM\cW(\widetilde{\cM})$
proves \eqref{2026092806}.

\subsection{Weighted matrix norms}\label{sec:weighted-norm-preliminaries}
For a matrix $\cM$ with rows in $\fX\subset\mathbb Z^d$ and columns
in $\fY\subset\mathbb Z^d$, define
\begin{equation}\label{2026092901}
 \matnorm{\cM}{\rho}
 =\sup_{y\in\fY}\sum_{x\in\fX}
       |\cM(x,y)|\langle x-y\rangle^\rho,
 \qquad \rho\ge0.
\end{equation}
If $\fY=\varnothing$, let $\matnorm{\cM}{\rho}=0$.
When $\matnorm{\cM}{0}<\infty$, $\cM$ defines a bounded map
from $\ell^1(\fY)$ to $\ell^1(\fX)$ with operator norm
$\matnorm{\cM}{0}$.

\begin{lemma}\label{lem:column-algebra}
Let $\rho\ge0$ and let $\fX,\fY,\fZ\subset\mathbb Z^d$.
For matrices $\cM:\fX\times\fY\to\mathbb C$ and
$\widetilde\cM:\fY\times\fZ\to\mathbb C$ with finite norms
in \eqref{2026092901},
\begin{equation}\label{2026092902}
 \matnorm{\cM\widetilde\cM}{\rho}
 \le\matnorm{\cM}{\rho}\matnorm{\widetilde\cM}{\rho}.
\end{equation}
For fixed $\fX$, matrices on $\fX\times\fX$ with finite
$\matnorm{\cdot}{\rho}$ form a Banach algebra.
If $\matnorm{\cM}{\rho}<1$ for a matrix on $\fX\times\fX$, then
\[
 (\Id_\fX-\cM)^{-1}=\sum_{k=0}^{\infty}\cM^k,
 \qquad
 \matnorm{(\Id_\fX-\cM)^{-1}}{\rho}
 \le\frac{1}{1-\matnorm{\cM}{\rho}}.
\]
\end{lemma}
\begin{proof}
The inequality
$\langle x-w\rangle\le\langle x-y\rangle\langle y-w\rangle$
implies
\begin{align*}
 \sum_{x\in\fX}|(\cM\widetilde\cM)(x,w)|\langle x-w\rangle^\rho
 &\le\sum_{y\in\fY}|\widetilde\cM(y,w)|\langle y-w\rangle^\rho
       \sum_{x\in\fX}|\cM(x,y)|\langle x-y\rangle^\rho\\
 &\le\matnorm{\cM}{\rho}\matnorm{\widetilde\cM}{\rho},
 \qquad w\in\fZ.
\end{align*}
Taking the supremum over $w$ proves \eqref{2026092902}.
If $(\cM_j)$ is Cauchy in $\matnorm{\cdot}{\rho}$, every entry has a limit
$\cM(x,y)$. By Fatou's lemma,
\[
 \matnorm{\cM_j-\cM}{\rho}
 \le\liminf_{k\to\infty}\matnorm{\cM_j-\cM_k}{\rho}\longrightarrow0.
\]
This proves completeness. If $\matnorm{\cM}{\rho}<1$, then
\eqref{2026092902} gives
$\sum_{k\ge0}\matnorm{\cM^k}{\rho}\le(1-\matnorm{\cM}{\rho})^{-1}$.
The Neumann series in this Banach algebra therefore gives the inverse
and the stated norm bound.
\end{proof}

\subsection{Finite-dimensional spectral facts}\label{sec:spectral-preliminaries}
For a Hermitian matrix $\cM$ on a finite set $\fF$ of cardinality
$m$, denote the eigenvalues, counted with multiplicity, by
\[
 \lambda_0(\cM)\le\cdots\le\lambda_{m-1}(\cM),
\]
and denote the spectrum by $\Spec\cM$.
We recall the min--max principle, Weyl's inequality and Cauchy
interlacing; see \cite{HornJohnson2013} and
\cite[Chapter~II]{Kato1995}.
For $0\le j<m$,
\begin{equation}\label{2026092903}
 \lambda_j(\cM)
 =\min_{\substack{V\subset\mathbb C^\fF\\\dim V=j+1}}
       \max_{\substack{v\in V\\\|v\|_2=1}}\langle v,\cM v\rangle.
\end{equation}
In particular, if $\cM,\widetilde{\cM}$ are Hermitian matrices on the same
set $\fF$ and $c\in\mathbb R$, then
\begin{align}
 \widetilde{\cM}\ge\cM+c\Id_\fF
 &\quad\Longrightarrow\quad
 \lambda_j(\widetilde{\cM})\ge\lambda_j(\cM)+c,
 \qquad 0\le j<m,\label{2026092904}\\
 |\lambda_j(\cM)-\lambda_j(\widetilde{\cM})|
 &\le\|\cM-\widetilde{\cM}\|_{\rm op},
 \qquad 0\le j<m.\label{2026092905}
\end{align}
The matrix order in \eqref{2026092904} means
$\langle v,(\widetilde{\cM}-\cM)v\rangle\ge c\|v\|_2^2$ for every
$v\in\mathbb C^\fF$.
If $\fP\subset\fF$ and $\#\fP=m-k\ge1$, then
\[
 \lambda_j(\cM)\le\lambda_j(\cM_\fP)\le\lambda_{j+k}(\cM),
 \qquad 0\le j<m-k.
\]
By \eqref{2026092905}, each ordered eigenvalue is continuous
under continuous matrix perturbations, including when eigenvalues coincide.

\section{The initial induction step}
\label{sec:initial-construction}
Assume $r>\tau$.
We first solve for the coordinates at which the original diagonal
entries are bounded away from zero. In a suitable finite region,
the remaining equation involves only $v(o)$.
We select a simple local eigenvalue by a homotopy in the coupling
constant and recover the corresponding eigenvector by back substitution.

\subsection{The equation and the first division of the lattice}\label{sec:initial-step}
Fix $\theta\in\fT$, $z\in\mathbb R$ and a finite set
$\fO\subset\mathbb Z^d$.  For $v\in\mathbb C^\fO$, the equation
$\cA_\fO(z;\theta)v=0$ reads
\begin{equation}\label{2026092906}
 \bigl(z-f(\theta+p\cdot\om)\bigr)v(p)
 -\eps\sum_{q\in\fO}\phi(p-q)v(q)=0,
 \qquad p\in\fO.
\end{equation}
Here $\cA_\fO$ is defined in \eqref{2026092802}.
The diagonal entries $z-f(\theta+p\cdot\om)$ can be arbitrarily small.
We first solve \eqref{2026092906} at points satisfying
$|z-f(\theta+p\cdot\om)|\ge\beta$, where we fix
\begin{equation}\label{2026093001}
 \beta=\frac{L\gamma}{32\cdot67200^\tau},\qquad
 \delta_0=\beta,\qquad E_0=f.
\end{equation}
Then we define
\begin{equation}
 \fS_0=\{p\in\mathbb Z^d:|z-f(\theta+p\cdot\om)|<\beta\},\qquad \fQ_0=\mathbb Z^d\setminus\fS_0.
 \label{2026093002}
\end{equation}
Let
\begin{equation}\label{2026093003}
 a=|\eps|\sup_{k\in\mathbb Z^d}\langle k\rangle^r|\phi(k)|.
\end{equation}
For $p,q\in\fO$ with $p\ne q$, assumption \eqref{2026092601} implies
\begin{equation}\label{2026093004}
 |\cA_\fO(z;\theta)(p,q)|
 =|\eps\phi(p-q)|\le a\langle p-q\rangle^{-r}.
\end{equation}

\subsection{Solving for the components in \texorpdfstring{$\fQ_0$}{Q0}}
Let $I_0=\fQ_0\cap\fO$ and $P_0=\fS_0\cap\fO$.
Write $\cA_\fO^{(-1)}=\cA_\fO(z;\theta)$.
Lemma~\ref{lem:initial-elimination} expresses $\left.v\right|_{I_0}$
in terms of $\left.v\right|_{P_0}$ and bounds
$\Sch_{I_0}(\cA_\fO^{(-1)})$.

Fix a constant $\mathrm C_{3}=\mathrm C_{3}(d,r)\ge16$ for the
sparse convolution and tail estimates in Lemma~\ref{lem:sparse-sums},
chosen also so that
\begin{equation}\label{2026093005}
 \sum_{k\in\mathbb Z^d}\langle k\rangle^{-r}\le\mathrm C_{3}.
\end{equation}
The proof of Lemma~\ref{lem:sparse-sums} gives a common finite bound
for the convolution, tail and full-lattice sums, depending only on $d,r$.

\begin{lemma}\label{lem:initial-elimination}
If $a\le(2\mathrm C_{3})^{-1}\beta$, then
$(\cA_\fO^{(-1)})_{I_0}$ is invertible for every
$\theta\in\fT$, $z\in\mathbb R$ and finite
$\fO\subset\mathbb Z^d$. Moreover, 
the equation \eqref{2026092906} is equivalent to
\begin{align}
 &\left.v\right|_{I_0}
 =- (\cA_\fO^{(-1)})_{I_0}^{-1}
       (\cA_\fO^{(-1)})_{I_0,P_0}\left.v\right|_{P_0},
 \label{2026093006}\\
 &\cA_\fO^{(0)}\left.v\right|_{P_0}=0,
 \label{2026100101}
\end{align}
where
\[
 \cA_\fO^{(0)}
 = (\cA_\fO^{(-1)})_{P_0}
 -(\cA_\fO^{(-1)})_{P_0,I_0}
       (\cA_\fO^{(-1)})_{I_0}^{-1}
       (\cA_\fO^{(-1)})_{I_0,P_0}.
\]
Furthermore, for $p\in P_0$,
\begin{equation}\label{2026100102}
 \left|\cA_\fO^{(0)}(p,p)
       -\bigl(z-f(\theta+p\cdot\om)\bigr)\right|
 \le2\mathrm C_{3}a^2\beta^{-1};
\end{equation}
and for distinct $p,q\in P_0$,
\begin{equation}\label{2026100103}
 |\cA_\fO^{(0)}(p,q)|
 \le a(1+2\mathrm C_{3} a\beta^{-1})\langle p-q\rangle^{-r}
 \le2a\langle p-q\rangle^{-r}.
\end{equation}
\end{lemma}
\begin{proof}
By \eqref{2026093002},
\begin{equation}\label{2026100104}
 |\cA_\fO^{(-1)}(p,p)|
 =|z-f(\theta+p\cdot\om)|\ge\beta,
 \qquad p\in I_0.
\end{equation}
Since $\fO\subset\mathbb Z^d$ is $1$-separated and
$a/\beta\le(2\mathrm C_{3})^{-1}$,
\eqref{2026093004} and \eqref{2026100104}
verify the assumptions of Lemma~\ref{lem:sparse-schur} with
$\fX=\fO$, $\fP=P_0$, $\fQ=I_0$, $M=1$,
$a_*=a$ and $\Delta=\beta$.
Thus $(\cA_\fO^{(-1)})_{I_0}$ is invertible.
The equations \eqref{2026092906} on $I_0$ imply
\eqref{2026093006}; substituting \eqref{2026093006}
into \eqref{2026092906} on $P_0$ proves
\eqref{2026100101}.
Estimates \eqref{2026102102} and \eqref{2026102101} prove
\eqref{2026100102} and the first inequality in
\eqref{2026100103}.  Since
$2\mathrm C_{3}a/\beta\le1$,
we obtain the second inequality in \eqref{2026100103}.
\end{proof}

\subsection{Construction of the first local eigenvalue and eigenvector}
Lemma~\ref{lem:initial-elimination} expresses $\left.v\right|_{I_0}$ in terms of
$\left.v\right|_{P_0}$, but \eqref{2026100101} remains to be solved
on $P_0$. We choose a fixed finite region around $o$ in which all
other components can be solved in terms of $v(o)$ for energies close
to $f(\theta)$. Let
\begin{equation}\label{2026100105}
 N_0=1,\qquad N_1=16,\qquad
 \fB_0=\{o\},\qquad \fB_1=\fL_{4N_1}.
\end{equation}
The set $\fL_{4N_1}$ is defined in \eqref{2026092705}.
For $p\in B_1^\circ$, periodicity and monotonicity in
\eqref{2026092603}--\eqref{2026092604},
together with \eqref{2026092502} and \eqref{2026093001}, imply
\begin{equation}\label{2026100201}
 |f(\theta+p\cdot\om)-f(\theta)|
 \ge L\|p\cdot\om\|_{\mathbb T}
 \ge L\gamma(1+4N_1)^{-\tau}>8\beta.
\end{equation}
Thus, with
\[
 \fJ_1(\theta)=[f(\theta)-4\beta,f(\theta)+4\beta],
\]
we obtain
\begin{equation}\label{2026100202}
 |z-f(\theta+p\cdot\om)|>4\beta,
 \qquad z\in\fJ_1(\theta),\quad p\in\fB_1^\circ.
\end{equation}
In particular, $\fB_1^\circ\subset\fQ_0$ for $z\in J_1(\theta)$.

To select the eigenvalue issuing from $f(\theta)$, we vary the
coupling constant. For $t\in[0,1]$, define
\begin{equation}\label{2026100203}
 \cH(\theta;t)
 =\operatorname{diag}_{p\in\mathbb Z^d}f(\theta+p\cdot\om)
                    +t\eps\cT.
\end{equation}
For a finite set $F\subset\mathbb Z^d$, let
\begin{equation}
 \cH_F(\theta;t)
 =\operatorname{diag}_{p\in F}f(\theta+p\cdot\om)+t\eps\cT_F,\qquad \cA_F(z;\theta;t)=z\Id_F-\cH_F(\theta;t).
 \label{2026100204}
\end{equation}
For fixed $\theta\in\Theta$, the matrix $\cH_F(\theta;t)$ is finite
and Hermitian, and
\[
 \cH_F(\theta;t)-\cH_F(\theta;s)
 =(t-s)\eps\cT_F,\qquad s,t\in[0,1].
\]
By \eqref{2026092905}, every ordered eigenvalue of
$\cH_F(\theta;t)$ is continuous in $t$.
At $t=1$, $\cH_F(\theta;1)=\cH_F(\theta)$ and
$\cA_F(z;\theta;1)=\cA_F(z;\theta)$, with the right-hand sides defined in
\eqref{2026092801} and \eqref{2026092802}.
The parameters $\beta,a$ and the set $B_1$ are fixed throughout
$t\in[0,1]$. We use $\cA_\Omega^{(h)}(z;\theta;t)$, $h=-1,0$,
for the matrices before and after the first elimination.
Since $|t\eps\phi(p-q)|\le a\langle p-q\rangle^{-r}$,
Lemma~\ref{lem:initial-elimination} yields
\eqref{2026093006}--\eqref{2026100103}
for every $t\in[0,1]$ whenever $a\le(2\mathrm C_{3})^{-1}\beta$.

Assume $a\le(2\mathrm C_{3})^{-1}\beta$ for the rest of this subsection.
By \eqref{2026100202} and
Lemma~\ref{lem:initial-elimination},
$\cA_{\fB_1^\circ}(z;\theta;t)$ is invertible for
$z\in J_1(\theta)$ and $t\in[0,1]$. Define
\begin{samepage}
\begin{equation}
 s_1(z;\theta;t)
 =z-f(\theta)
 -c_1(t)^*\cA_{\fB_1^\circ}(z;\theta;t)^{-1}c_1(t),\qquad c_1(t)=-t\eps\Pi_{\fB_1^\circ}\cT e_o.
 \label{2026100205}
\end{equation}
Here $s_1$ is the Schur complement in
\eqref{2026092803} and $c_1(t)$ is the column of
$\cA_{B_1}(z;\theta;t)$ with row set $B_1^\circ$
and column index $o$. We write $s_1(z;\theta)=s_1(z;\theta;1)$.
\end{samepage}
For $z\in J_1(\theta)$, the equation
$\cA_{B_1}(z;\theta;t)v=0$ is equivalent to
\begin{equation*}
 \left.v\right|_{\fB_1^\circ}
 =-\cA_{\fB_1^\circ}(z;\theta;t)^{-1}c_1(t)v(o),\qquad s_1(z;\theta;t)v(o)=0.
\end{equation*}
Every nonzero solution of $\cA_{B_1}(z;\theta;t)v=0$ has $v(o)\ne0$. Hence
\[
 z\in\Spec\cH_{B_1}(\theta;t)
 \quad\Longleftrightarrow\quad s_1(z;\theta;t)=0,
 \qquad z\in J_1(\theta).
\]
Block multiplication gives
\begin{equation}\label{2026100206}
 \det\cA_{B_1}(z;\theta;t)
 =\det\cA_{B_1^\circ}(z;\theta;t)\,s_1(z;\theta;t).
\end{equation}
Since $a/\beta\le(2\mathrm C_{3})^{-1}$,
the bounds \eqref{2026093004} and
\eqref{2026100202} verify the assumptions of
Lemma~\ref{lem:sparse-schur} for $\cA_{B_1}(z;\theta;t)$ with
$\fX=B_1$, $\fP=\{o\}$, $\fQ=B_1^\circ$,
$M=1$, $a_*=a$ and $\Delta=\beta$.
Thus \eqref{2026102102} implies
\begin{equation}\label{2026100301}
 |s_1(z;\theta;t)-(z-f(\theta))|
 \le2\mathrm C_{3}a^2\beta^{-1},
 \qquad z\in J_1(\theta),\quad t\in[0,1].
\end{equation}
\begin{proposition}\label{prop:initial-step}
Assume $a\le(2\mathrm C_{3})^{-1}\beta$.
For $\theta\in\Theta$, let
$j=\#\{p\in B_1:f(\theta+p\cdot\om)<f(\theta)\}$ and define
\begin{equation}\label{2026100302}
 E_1(\theta;t)=\lambda_j(\cH_{B_1}(\theta;t)),
 \qquad t\in[0,1],\qquad E_1(\theta)=E_1(\theta;1).
\end{equation}
Then $E_1(\cdot;t)$ is real and $1$-periodic,
$t\mapsto E_1(\theta;t)$ is continuous and
$E_1(\theta;0)=f(\theta)$.
For every $\theta\in\Theta$ and $t\in[0,1]$, $E_1(\theta;t)$ is the unique
zero of $s_1(\cdot;\theta;t)$ in $J_1(\theta)$ and satisfies
\begin{equation}\label{2026100303}
 |E_1(\theta;t)-f(\theta)|
 \le2\mathrm C_{3}a^2\beta^{-1}\le\beta/16.
\end{equation}
The eigenvalue $E_1(\theta;t)$ is simple, with
\begin{equation}
 \Spec\cH_{B_1}(\theta;t)\cap J_1(\theta)=\{E_1(\theta;t)\},\qquad \dist\bigl(E_1(\theta;t),\Spec\cH_{B_1}(\theta;t)
                    \setminus\{E_1(\theta;t)\}\bigr)\ge3\beta.
 \label{2026100304}
\end{equation}
\end{proposition}
\begin{proof}
Fix $\theta\in\Theta$. By \eqref{2026100301} and
\eqref{2026093005},
\[
 s_1(f(\theta)-4\beta;\theta;t)<0
 <s_1(f(\theta)+4\beta;\theta;t),\qquad t\in[0,1].
\]
Since $\cA_{B_1^\circ}(z;\theta;t)$ is invertible on
$J_1(\theta)$, \eqref{2026100206} implies
\[
 f(\theta)\pm4\beta\notin\Spec\cH_{B_1}(\theta;t),
 \qquad t\in[0,1].
\]
By \eqref{2026092905} and
$f(\theta)\pm4\beta\notin\Spec\cH_{B_1}(\theta;t)$,
the functions $t\mapsto\lambda_k(\cH_{B_1}(\theta;t))-f(\theta)\pm4\beta$
are continuous and nonzero on $[0,1]$, hence have constant sign,
for every $0\le k<\#B_1$.
At $t=0$, \eqref{2026100201} places exactly the
eigenvalue with index $j$ in $J_1(\theta)$. Hence
\[
 \{k:\lambda_k(\cH_{B_1}(\theta;t))\in J_1(\theta)\}
 =\{j\},\qquad t\in[0,1].
\]
The eigenvalues are counted with multiplicity.
Consequently, $E_1(\theta;t)$ is simple and is the unique eigenvalue
of $\cH_{B_1}(\theta;t)$ in $J_1(\theta)$.
By \eqref{2026100206}, $s_1(\cdot;\theta;t)$ has the unique
zero $E_1(\theta;t)$ in $J_1(\theta)$.
Taking $z=E_1(\theta;t)$ in \eqref{2026100301} proves
\eqref{2026100303}. Since
$\Spec\cH_{B_1}(\theta;t)\setminus\{E_1(\theta;t)\}$ lies outside
$J_1(\theta)$,
\[
 \dist\bigl(E_1(\theta;t),\Spec\cH_{B_1}(\theta;t)
                  \setminus\{E_1(\theta;t)\}\bigr)
 \ge4\beta-2\mathrm C_{3}a^2\beta^{-1}\ge3\beta.
\]
This proves \eqref{2026100304}.
Since $\cH_{B_1}(\theta;t)$ and
$j=\#\{p\in B_1:f(\theta+p\cdot\om)<f(\theta)\}$ are $1$-periodic
in $\theta$, \eqref{2026100302} proves that
$E_1(\cdot;t)$ is $1$-periodic.
\end{proof}

Define the zero-extended vector $v_{1,\theta}(t)$ by
\begin{equation}
 (v_{1,\theta}(t))(o)=1,\qquad \left.v_{1,\theta}(t)\right|_{B_1^c}=0,\qquad \left.v_{1,\theta}(t)\right|_{B_1^\circ}
 =-\cA_{B_1^\circ}(E_1(\theta;t);\theta;t)^{-1}c_1(t).
 \label{2026100305}
\end{equation}
By \eqref{2026100205} and
$s_1(E_1(\theta;t);\theta;t)=0$,
\[
 \cH_{B_1}(\theta;t)\left.v_{1,\theta}(t)\right|_{B_1}
 =E_1(\theta;t)\left.v_{1,\theta}(t)\right|_{B_1}.
\]
Formula \eqref{2026100305} also proves continuity of
$t\mapsto v_{1,\theta}(t)$ and $v_{1,\theta}(0)=e_o$. Write $v_{1,\theta}=v_{1,\theta}(1)$
and $v_{0,\theta}(t)=e_o$.
Fix $0<\rho<r-\tau$ and define
\begin{equation}\label{2026100306}
 \|u\|_{1,\rho}=\sum_{m\in\mathbb Z^d}|u(m)|\langle m\rangle^\rho.
\end{equation}
For each $0\le s<r-d$, fix a finite coefficient
$\mathrm C_{4}(d,r,s)>0$ for \eqref{2026102104} with $\rho=s$.
For the chosen $\rho$, set
\begin{equation}\label{2026100307}
 \mathrm C_{4}=\mathrm C_{4}(d,r,\rho).
\end{equation}
For the decomposition $B_1=\{o\}\sqcup B_1^\circ$,
the assumptions of Lemma~\ref{lem:sparse-schur} were verified before
\eqref{2026100301}. Since $\rho<r-\tau<r-d$,
\eqref{2026102104} implies
\begin{equation}\label{2026100401}
 \sup_{\theta\in\Theta,\,t\in[0,1]}
 \|v_{1,\theta}(t)-v_{0,\theta}(t)\|_{1,\rho}
 \le\frac{\mathrm C_{4}}{4}\,a\beta^{-1}.
\end{equation}
We call $E_1$ the first Rellich function. To select the next coordinates
to solve in \eqref{2026100101} using
$E_1(\theta+p\cdot\om)$, we prove that $E_1$ is Lipschitz monotone
in the phase.

\subsection{Reordered eigenvalue branches}\label{sec:reordered-branches}
Let $\fF\subset\mathbb Z^d$ be finite with $o\in\fF$ and $m=\#\fF$.
For $t\in[0,1]$, order the eigenvalues of
$\cH_\fF(\theta;t)$ increasingly, counting multiplicities,
as in Subsection~\ref{sec:spectral-preliminaries}.
Arrange the phases $\theta\in[0,1)$ for which
$\theta+p\cdot\om\in\mathbb Z$ for some $p\in\fF$ as
\begin{equation}
 0=\vartheta_0<\vartheta_1<\cdots<\vartheta_{m-1}
                    <\vartheta_m=1,\qquad \{\vartheta_0,\ldots,\vartheta_{m-1}\}
 =\{\{-p\cdot\om\}:p\in\fF\}.
 \label{2026100402}
\end{equation}
For the fractional parts defined in \eqref{2026092702},
\eqref{2026092502} implies that $\{-p\cdot\om\}$, $p\in\fF$, are pairwise
distinct, and $o\in\fF$ implies $\vartheta_0=0$.
For $0\le j<m$, the functions $\theta\mapsto f(\theta+p\cdot\om)$,
$p\in\fF$, are Lipschitz monotone with constant $L$ on
$(\vartheta_j,\vartheta_{j+1})$.
By \eqref{2026092903}, each function
$\theta\mapsto\lambda_k(\cH_F(\theta;t))$, $0\le k<m$,
is Lipschitz monotone with constant $L$ on $(\vartheta_j,\vartheta_{j+1})$.
At $\vartheta_j=\{-p\cdot\om\}$, $1\le j<m$, the fractional part
$\{\theta+p\cdot\om\}$ returns to $0$.
We therefore increase the selected eigenvalue index by one at each
$\vartheta_j$, $1\le j<m$:
\begin{equation}\label{2026100403}
 \lambda(\fF;\theta;t)=\lambda_j(\cH_\fF(\theta;t)),
 \qquad
 \theta\in[\vartheta_j,\vartheta_{j+1})\cap\fT,
 \quad 0\le j<m.
\end{equation}
Extend $\lambda(\fF;\cdot;t)$ $1$-periodically to $\fT$ and write
$\lambda(\fF;\theta)=\lambda(\fF;\theta;1)$.
If $\theta\in[\vartheta_j,\vartheta_{j+1})\cap\Theta$, then exactly
$j$ diagonal entries of $\cH_F(\theta;0)$ are smaller than $f(\theta)$:
for $p\ne o$,
\[
 \{\theta+p\cdot\om\}<\theta
 \quad\Longleftrightarrow\quad \{-p\cdot\om\}\le\theta.
\]
Thus the eigenvalue branch associated with the diagonal entry $f(\theta)$
in \cite[Remark~4.2]{KPS2025} is precisely
\eqref{2026100403}. Hence
\begin{equation}\label{2026100404}
 \lambda(F;\theta;0)=f(\theta),\qquad \theta\in\Theta.
\end{equation}

\begin{proposition}[{\cite[Proposition~4.1]{KPS2025}}]\label{prop:kps-label}
For every finite $\fF\subset\mathbb Z^d$ containing $o$ and every
$t\in[0,1]$, the function $\lambda(\fF;\cdot;t)$ in
\eqref{2026100403} satisfies
\begin{equation}\label{2026100501}
 \lambda(\fF;\theta_2;t)-\lambda(\fF;\theta_1;t)\ge L(\theta_2-\theta_1),
 \qquad \theta_1,\theta_2\in\fT\cap[0,1),\quad \theta_1<\theta_2.
\end{equation}
\end{proposition}
For fixed $t$, the hopping matrix $t\eps\cT_F$ is Hermitian
and independent of $\theta$, as required in
\cite[Remark~4.2]{KPS2025}.
Since $f/L$ is Lipschitz monotone with constant $1$,
\cite[Proposition~4.1]{KPS2025} applied to $\cH_F(\theta;t)/L$
proves \eqref{2026100501}.

Taking $F=B_1$ in \eqref{2026100403} and using
\eqref{2026100302}, we obtain
\begin{equation}\label{2026100502}
 E_1(\theta;t)=\lambda(B_1;\theta;t),
 \qquad \theta\in\Theta,\quad t\in[0,1].
\end{equation}
Proposition~\ref{prop:kps-label} thus proves
$E_1(\theta_2;t)-E_1(\theta_1;t)\ge L(\theta_2-\theta_1)$ for
$\theta_1,\theta_2\in\Theta\cap[0,1)$ with $\theta_1<\theta_2$.

\subsection{Homotopy on a moving energy interval}
\label{sec:homotopy-labels}
Proposition~\ref{prop:initial-step} counts the eigenvalues of
$\cH_{B_1}(\theta;t)$ in the fixed interval $J_1(\theta)$.
For the next local eigenvalue, the energy interval is centered at
$E_1(\theta;t)$ and therefore varies with $t$.
Proposition~\ref{prop:homotopy-monotonicity} extends the eigenvalue count
to an interval with a continuous center $h(\theta;t)$.

\begin{proposition}\label{prop:homotopy-monotonicity}
Let $F\subset\mathbb Z^d$ be finite with $o\in F$ and let
$\eta>0$. For each $\theta\in\Theta$, let
$t\mapsto h(\theta;t)$ be continuous on $[0,1]$, with
$h(\theta;0)=f(\theta)$. Suppose that
\begin{align*}
 &|f(\theta+p\cdot\om)-f(\theta)|>\eta,\quad p\in F^\circ,\quad\theta\in\Theta,\\
 &h(\theta;t)\pm\eta\notin\Spec\cH_F(\theta;t),\quad \theta\in\Theta,\quad t\in[0,1].
\end{align*}
Then $[h(\theta;t)-\eta,h(\theta;t)+\eta]$ contains exactly one
eigenvalue of $\cH_F(\theta;t)$, counted with multiplicity, namely
$\lambda(F;\theta;t)$ from \eqref{2026100403}.
The function $\lambda(F;\theta;t)$ is continuous in $t$, is $1$-periodic
in $\theta$ and satisfies \eqref{2026100501}.
\end{proposition}
\begin{proof}
Fix $\theta\in\Theta$ and let
$j=\#\{p\in F:f(\theta+p\cdot\om)<f(\theta)\}$.
For each $0\le k<\#F$, the functions
$t\mapsto\lambda_k(\cH_F(\theta;t))-h(\theta;t)\pm\eta$
are continuous and nonzero on $[0,1]$.
The sign argument in the proof of Proposition~\ref{prop:initial-step},
together with the eigenvalue count at $t=0$, gives
\[
 \{k:|\lambda_k(\cH_F(\theta;t))-h(\theta;t)|<\eta\}
 =\{j\},\qquad t\in[0,1].
\]
The definition \eqref{2026100403} and the count
$j=\#\{p\in F:f(\theta+p\cdot\om)<f(\theta)\}$ imply
$\lambda_j(\cH_F(\theta;t))=\lambda(F;\theta;t)$.
By \eqref{2026092905}, $t\mapsto\lambda(F;\theta;t)$ is continuous.
The periodic extension in \eqref{2026100403} makes
$\lambda(F;\cdot;t)$ $1$-periodic, and
Proposition~\ref{prop:kps-label} proves \eqref{2026100501}.
\end{proof}

To verify $h(\theta;t)\pm\eta\notin\Spec\cH_F(\theta;t)$ in
Proposition~\ref{prop:homotopy-monotonicity}, we compare the scalar
Schur complement at two energies.

\begin{lemma}\label{lem:schur-energy-comparison}
Let $F\subset\mathbb Z^d$ be finite with $o\in F$ and fix
$\theta\in\Theta$ and $t\in[0,1]$. Suppose that
$\cA_{F^\circ}(z;\theta;t)$ is invertible for every
$z$ in a real interval $J$. Denote the scalar Schur complement by
\[
 s(z;\theta;t)=\Sch_{F^\circ}(\cA_F(z;\theta;t)).
\]
For $z_1<z_2$ in $J$,
\begin{equation}\label{2026100503}
 s(z_2;\theta;t)-s(z_1;\theta;t)\ge z_2-z_1.
\end{equation}
If $E\in J$ and $s(E;\theta;t)=0$, then
$|s(z;\theta;t)|\ge|z-E|$ for every $z\in J$.
\end{lemma}
\begin{proof}
Let $c=-t\eps\Pi_{F^\circ}\cT e_o$ and omit the fixed
arguments $\theta,t$. Since
$\Spec\cH_{F^\circ}(\theta;t)\cap J=\varnothing$, for every
$\mu\in\Spec\cH_{F^\circ}(\theta;t)$ we have
$((z_2-\mu)(z_1-\mu))^{-1}>0$.
Thus $\cA_{F^\circ}(z_2)^{-1}\cA_{F^\circ}(z_1)^{-1}$ is positive
semidefinite, and the resolvent identity gives
\begin{equation*}
 s(z_2)-s(z_1)
 =(z_2-z_1)\bigl(1+c^*\cA_{F^\circ}(z_2)^{-1}
                              \cA_{F^\circ}(z_1)^{-1}c\bigr) \ge z_2-z_1.
\end{equation*}
Taking $z_1=E$ when $z>E$ and $z_2=E$ when $z<E$ in
\eqref{2026100503} proves
$|s(z;\theta;t)|\ge|z-E|$.
\end{proof}

Since Proposition~\ref{prop:initial-step} identifies
$E_1(\theta;t)$ as the zero of $s_1(\cdot;\theta;t)$ in $J_1(\theta)$,
Lemma~\ref{lem:schur-energy-comparison} with $F=B_1$ implies
\begin{equation}\label{2026100504}
 |s_1(z;\theta;t)|\ge|z-E_1(\theta;t)|,
 \qquad z\in J_1(\theta),\quad t\in[0,1].
\end{equation}

\subsection{Continuing the elimination}
\label{sec:next-equations}
Fix $\theta\in\Theta$, $z\in\mathbb R$, $t\in[0,1]$ and a finite
set $\Omega\subset\mathbb Z^d$, and assume $a\le(2\mathrm C_{3})^{-1}\beta$.
For fixed $t$, we omit $t$ from the matrices, Schur complements
and eigenvalue functions.
By Lemma~\ref{lem:initial-elimination}, the remaining equation
\eqref{2026100101} is
\[
 \cA_\Omega^{(0)}\left.v\right|_{P_0}=0,
 \qquad P_0=S_0\cap\Omega.
\]
The diagonal entries $\cA_\Omega^{(0)}(p,p)$ depend on $\Omega$.
We first examine the equation in $B_1(p)$, defined in
\eqref{2026092704}.

For $p\in S_0$, we have $z\in J_1(\theta+p\cdot\om)$, so
\eqref{2026100202} implies
$B_1(p)\cap S_0=\{p\}$.
By \eqref{2026100205}, \eqref{2026092708}
and \eqref{2026100504},
\begin{equation*}
 \cA_{B_1(p)}^{(0)}(p,p)=s_1(z;\theta+p\cdot\om),\qquad |s_1(z;\theta+p\cdot\om)|\ge|z-E_1(\theta+p\cdot\om)|.
\end{equation*}
We therefore select the next coordinates using the distance
to $E_1(\theta+p\cdot\om)$. With $N_1=16$ as in
\eqref{2026100105}, define
\begin{align}
 &\delta_1=\beta N_1^{-\tau},\label{2026100505}\\
 &S_1=\{p\in S_0:|z-E_1(\theta+p\cdot\om)|<\delta_1\},
 \qquad Q_1=S_0\setminus S_1.\label{2026100506}
\end{align}
In particular, $|s_1(z;\theta+p\cdot\om)|\ge\delta_1$ for
$p\in Q_1$. The points in $S_1$ satisfy a separation estimate
that follows from the next lemma.

\begin{lemma}\label{lem:branch-separation}
Let $E:\Theta\to\mathbb R$ be $1$-periodic and satisfy
$E(\theta_2)-E(\theta_1)\ge L(\theta_2-\theta_1)$ for $\theta_1,\theta_2\in\Theta\cap[0,1)$ with $\theta_1<\theta_2$.
For distinct $p,q\in\mathbb Z^d$ and $\theta\in\Theta$,
\begin{equation}\label{2026100507}
 |E(\theta+p\cdot\om)-E(\theta+q\cdot\om)|
 \ge L\|(p-q)\cdot\om\|_{\mathbb T}
 \ge L\gamma\langle p-q\rangle^{-\tau}.
\end{equation}
\end{lemma}
\begin{proof}
Using the periodicity identity in \eqref{2026092603}
and the Lipschitz monotonicity inequality \eqref{2026092604},
both assumed here for $E$, we obtain
\begin{equation*}
 |E(\theta+p\cdot\om)-E(\theta+q\cdot\om)| \ge L|\{\theta+p\cdot\om\}-\{\theta+q\cdot\om\}|
 \ge L\|(p-q)\cdot\om\|_{\mathbb T}.
\end{equation*}
The Diophantine condition \eqref{2026092502} now proves \eqref{2026100507}.
\end{proof}
By \eqref{2026100502},
Proposition~\ref{prop:kps-label} and Lemma~\ref{lem:branch-separation},
distinct $p,q\in S_1$ satisfy
\begin{equation}\label{2026100508}
 \langle p-q\rangle>
 \left(\frac{L\gamma}{2\beta}\right)^{1/\tau}N_1.
\end{equation}

Let $I_1=Q_1\cap\Omega$ and $P_1=S_1\cap\Omega$.
The rows indexed by $I_1$ in \eqref{2026100101} read
\begin{equation}\label{2026100509}
 (\cA_\Omega^{(0)})_{I_1}\left.v\right|_{I_1}
 +(\cA_\Omega^{(0)})_{I_1,P_1}\left.v\right|_{P_1}=0.
\end{equation}
To solve \eqref{2026100509} for $\left.v\right|_{I_1}$, we need
$(\cA_\Omega^{(0)})_{I_1}$ to be invertible.
The off-diagonal entries are already bounded by
\eqref{2026100103}.
The local lower bound $|s_1(z;\theta+p\cdot\om)|\ge\delta_1$
will imply $|\cA_\Omega^{(0)}(p,p)|\ge\delta_1/2$ if we prove
\begin{equation}\label{2026100510}
 |\cA_\Omega^{(0)}(p,p)-s_1(z;\theta+p\cdot\om)|
 \le\delta_1/2,\qquad p\in I_1.
\end{equation}
Under \eqref{2026100510} and $4a/\delta_1\le(2\mathrm C_{3})^{-1}$,
we verify the assumptions of Lemma~\ref{lem:sparse-schur} with
$\fX=P_0$, $\fP=P_1$, $\fQ=I_1$, $M=1$, $a_*=2a$ and
$\Delta=\delta_1/2$.
Thus $(\cA_\Omega^{(0)})_{I_1}$ is invertible, and
\eqref{2026100101} becomes
\begin{equation}
 \left.v\right|_{I_1}=-(\cA_\Omega^{(0)})_{I_1}^{-1}
       (\cA_\Omega^{(0)})_{I_1,P_1}\left.v\right|_{P_1},\qquad \cA_\Omega^{(1)}\left.v\right|_{P_1}=0,\qquad
 \cA_\Omega^{(1)}=\Sch_{I_1}(\cA_\Omega^{(0)}).
 \label{2026100601}
\end{equation}

To compare $\cA_\Omega^{(0)}(p,p)$ with
$\cA_{B_1(p)}^{(0)}(p,p)$, we require $B_1(p)\subset\Omega$.
Since $B_1(p)\cap S_0=\{p\}$, the eliminated coordinates then satisfy
\[
 Q_0\cap\Omega=B_1^\circ(p)
 \sqcup\bigl(Q_0\cap(\Omega\setminus B_1(p))\bigr).
\]
Every point of $Q_0\cap(\Omega\setminus B_1(p))$ has distance
greater than $4N_1$ from $p$.
The distance bound and \eqref{2026093004} will be used
to estimate the effect of the additional eliminated coordinates.
The comparison is also needed for the remaining coordinates:
by \eqref{2026100601}, for $p\in P_1$,
\[
 \cA_\Omega^{(1)}(p,p)=\cA_\Omega^{(0)}(p,p)
 -(\cA_\Omega^{(0)})_{p,I_1}
   (\cA_\Omega^{(0)})_{I_1}^{-1}(\cA_\Omega^{(0)})_{I_1,p}.
\]
At $z=E_1(\theta+p\cdot\om)$, the identity
$s_1(z;\theta+p\cdot\om)=0$ allows the region comparison
to bound $|\cA_\Omega^{(0)}(p,p)|$ in the remaining equation.
We therefore impose
\begin{equation}\label{2026100602}
 B_1(p)\subset\Omega,\qquad p\in S_0(z;\theta)\cap\Omega.
\end{equation}
On a larger finite region, we seek a local eigenvalue near
$E_1(\theta;t)$ by solving for all coordinates except $v(o)$.
For the next application of Proposition~\ref{prop:homotopy-monotonicity}
with center $E_1(\theta;t)$, the finite region must be fixed
as $t$ varies.
We arrange \eqref{2026100602} uniformly for
$\theta\in\Theta$, $t\in[0,1]$ and $|z-E_1(\theta;t)|\le4\delta_1$.
For $p\in S_0(z;\theta)$, \eqref{2026100303} gives
\begin{equation*}
 L\|p\cdot\om\|_{\mathbb T}
 \le |f(\theta+p\cdot\om)-f(\theta)| <\beta+4\delta_1+\beta/16<8\beta.
\end{equation*}
Define
\begin{equation}\label{2026100603}
 Z_0=\{p\in\mathbb Z^d:L\|p\cdot\om\|_{\mathbb T}\le8\beta\}.
\end{equation}
Then
\begin{equation}\label{2026100604}
 S_0(z;\theta)\subset Z_0,
 \qquad |z-E_1(\theta;t)|\le4\delta_1,
\end{equation}
for every $\theta\in\Theta$ and $t\in[0,1]$.
The set $Z_0$ is independent of $z$, $\theta$ and $t$.
By \eqref{2026100604}, a fixed finite set $F\supset B_1$ satisfies
\eqref{2026100602} for every $\theta\in\Theta$,
$t\in[0,1]$ and $|z-E_1(\theta;t)|\le4\delta_1$, provided that
\[
 p\in F\cap Z_0\quad\Longrightarrow\quad B_1(p)\subset F.
\]
We formulate the resulting induction hypotheses in
Section~\ref{sec:inductive-hypotheses}.

\section{Inductive hypotheses}\label{sec:finite-volume}
\label{sec:inductive-hypotheses}
At step $n\ge1$, extending the construction in
Section~\ref{sec:initial-construction},
the eigenvalues $E_m$ and corresponding eigenvectors $v_{m,\theta}$
are known for $0\le m\le n$, and the coordinates in
$Q_0,\ldots,Q_{n-1}$ have been solved.
The next step solves for the coordinates in $Q_n$, constructs
$E_{n+1}$ and the corresponding eigenvector $v_{n+1,\theta}$,
and estimates
$E_{n+1}-E_n$ and $v_{n+1,\theta}-v_{n,\theta}$.

\subsection{Scales and the equations at step \texorpdfstring{$n$}{n}}\label{sec:multiscale-notation}
\emph{Scales and fixed sets.}
Retain $\beta$ from \eqref{2026093001} and $a$ from
\eqref{2026093003}. Extend \eqref{2026100105}
and \eqref{2026100505} by
\begin{align}
 N_h&=16^h,\qquad \delta_h=\beta N_h^{-\tau},\qquad h\ge0,
 \label{2026100605}\\
 M_0&=1,\qquad M_h=N_{h-1}\quad(h\ge1),\qquad
 \eta_h=a\delta_h^{-1}M_h^{-r}\quad(h\ge0).
 \label{2026100606}
\end{align}
The numbers $N_h$ specify the sizes of the finite regions, $\delta_h$
specify the energy thresholds and $M_h$ specify the separation used
in the matrix estimates.  Since $r>\tau$, we have
\begin{equation}\label{2026100701}
 \sum_{h\ge0}\eta_h
 =a\beta^{-1}\left(1+16^\tau
          \sum_{h\ge1}N_{h-1}^{\tau-r}\right)<\infty.
\end{equation}

Extend the set $Z_0$ in \eqref{2026100603} by defining
\begin{equation}\label{2026100702}
 Z_h=\{p\in\mathbb Z^d:L\|p\cdot\om\|_{\mathbb T}\le8\delta_h\},
 \qquad h\ge0.
\end{equation}
For every $h\ge0$, the set $Z_h$ is independent of $z$, $\theta$
and $t$, and satisfies $Z_{h+1}\subset Z_h$ and $o\in Z_h$. For distinct $p,q\in Z_h$, the Diophantine condition \eqref{2026092502}
and $\langle p-q\rangle\le2|p-q|$ imply
\begin{equation}\label{2026100703}
 |p-q|\ge\frac12\left(\frac{L\gamma}{16\beta}\right)^{1/\tau}N_h
 >33600N_h.
\end{equation}
Here $\frac12(L\gamma/(16\beta))^{1/\tau}>33600$ follows from
the fixed choice \eqref{2026093001}.

\emph{The remaining equations.}
We first write the notation for $t=1$.
Given $E_0=f,E_1,\ldots,E_n$ and finite sets $B_0,\ldots,B_n$
with $B_0=\{o\}$, fix $\theta\in\Theta$ and $z\in\mathbb R$.
For $0\le h\le n$, define
\begin{align}
 \Res_h&=\{p\in\mathbb Z^d:|E_h(\theta+p\cdot\om)-z|<\delta_h\},
 \label{2026100704}\\
 S_{-1}&=\mathbb Z^d,\qquad S_h=\bigcap_{\ell=0}^{h}\Res_\ell,
 \label{2026100705}\\
 Q_h&=S_{h-1}\setminus S_h.
 \label{2026100706}
\end{align}
For $h=0,1$, \eqref{2026100705} and
\eqref{2026100706} agree with
\eqref{2026093002} and \eqref{2026100506}.
Whenever $E_h$ satisfies the hypotheses of
Lemma~\ref{lem:branch-separation}, distinct $p,q\in S_h$ satisfy
\[
 L\|(p-q)\cdot\om\|_{\mathbb T}
 \le |E_h(\theta+p\cdot\om)-E_h(\theta+q\cdot\om)|
 <2\delta_h,
 \qquad q-p\in Z_h.
\]
For a finite set $\Omega$, let
\begin{equation}\label{2026100707}
 I_h=Q_h\cap\Omega,\qquad P_h=S_h\cap\Omega,\qquad P_{-1}=\Omega.
\end{equation}
For $0\le h\le n$, define $\cA_\Omega^{(h)}$ recursively by
\begin{equation}\label{2026100708}
 \cA_\Omega^{(-1)}=\cA_\Omega(z;\theta),\qquad
 \cA_\Omega^{(h)}=\Sch_{I_h}(\cA_\Omega^{(h-1)}),
\end{equation}
provided that $(\cA_\Omega^{(h-1)})_{I_h}$ is invertible.
The matrix $\cA_\Omega^{(h)}$ has row and column set $P_h$.
The Schur identity \eqref{2026092803} makes
$\cA_\Omega^{(h-1)}\left.v\right|_{P_{h-1}}=0$ equivalent to
\begin{align}
 &\left.v\right|_{I_h}=-(\cA_\Omega^{(h-1)})_{I_h}^{-1}
             (\cA_\Omega^{(h-1)})_{I_h,P_h}\left.v\right|_{P_h},
 \label{2026100801}\\
 &\cA_\Omega^{(h)}\left.v\right|_{P_h}=0.
 \label{2026100802}
\end{align}
After eliminating $I_0,\ldots,I_{n-1}$, the remaining equation is
$\cA_\Omega^{(n-1)}\left.v\right|_{P_{n-1}}=0$.
Since $P_{n-1}=I_n\sqcup P_n$, the next elimination requires
$(\cA_\Omega^{(n-1)})_{I_n}$ to be invertible.

\emph{The finite regions.}
For $1\le i\le n$, comparing $\cA_{B_i(p)}^{(i-1)}(p,p)$ with
$\cA_\Omega^{(i-1)}(p,p)$ requires $B_i(p)\subset\Omega$
for $p\in S_{i-1}\cap\Omega$, as in
\eqref{2026100602} for $i=1$.
To construct $E_{n+1}$, the preceding eliminations must also be
available on the new sets $B_{n+1}$ and $B_{n+1}^\circ$.
We therefore formulate the estimates for every finite $\Omega$
satisfying the required inclusions for the already constructed
$B_1,\ldots,B_i$; the region $\Omega$ can then be chosen after
the estimates have been proved.
Let
\begin{equation}\label{2026100803}
 \Upsilon_0=\{\Omega\subset\mathbb Z^d:\#\Omega<\infty\}.
\end{equation}
For $1\le i\le n$, define recursively
\[
 \Upsilon_i=\{\Omega\in\Upsilon_{i-1}:
 B_i(p)\subset\Omega\text{ for every }p\in S_{i-1}\cap\Omega\}.
\]
Equivalently, $\Omega\in\Upsilon_i$ if and only if $\Omega$ is finite and
\begin{equation}\label{2026100804}
 B_h(p)\subset\Omega,
 \qquad 1\le h\le i,\quad p\in S_{h-1}\cap\Omega.
\end{equation}
For $i=1$, condition \eqref{2026100804} reduces to
\eqref{2026100602}.
The condition $\Omega\in\Upsilon_i$ does not require
$S_i\cap\Omega=\varnothing$; the coordinates in $S_i\cap\Omega$
remain in \eqref{2026100802}.
The sets $\Res_h,S_h,Q_h$ are defined on $\mathbb Z^d$ before
$\Omega$ is chosen. Through the sets $S_{h-1}(z;\theta)$ in
\eqref{2026100804}, the family
$\Upsilon_i=\Upsilon_i(z;\theta)$ depends on $z$ and $\theta$.

\emph{The local Schur complement.}
For $1\le m\le n$, let $c_m=-\eps\Pi_{B_m^\circ}\cT e_o$.
For $z$ such that $\cA_{B_m^\circ}(z;\theta)$ is invertible, define
\begin{equation}\label{2026100901}
 s_m(z;\theta)=z-f(\theta)
       -c_m^*\cA_{B_m^\circ}(z;\theta)^{-1}c_m.
\end{equation}
For $m=1$, definition \eqref{2026100901} agrees with
\eqref{2026100205}.
As in \eqref{2026100206},
\begin{equation}\label{2026100902}
 \det\cA_{B_m}(z;\theta)
 =\det\cA_{B_m^\circ}(z;\theta)\,s_m(z;\theta).
\end{equation}
For $1\le m\le n$, let
\begin{equation}\label{2026100903}
 J_m(\theta)=[E_{m-1}(\theta)-4\delta_{m-1},
                  E_{m-1}(\theta)+4\delta_{m-1}].
\end{equation}

\emph{The homotopy parameter.}
For $t\in[0,1]$, use $\cA_F(z;\theta;t)$ from
\eqref{2026100204}, $E_h(\theta;t)$ and
$c_m(t)=-t\eps\Pi_{B_m^\circ}\cT e_o$ in
\eqref{2026100704}--\eqref{2026100903}.
We write the resulting functions, matrices and sets as
$s_m(z;\theta;t)$, $J_m(\theta;t)$,
$\cA_\Omega^{(h)}(z;\theta;t)$, $\Res_h(z;\theta;t)$,
$S_h(z;\theta;t)$, $Q_h(z;\theta;t)$ and
$\Upsilon_i(z;\theta;t)$.
The arguments $z,\theta,t$ are omitted when fixed.
The sets $B_h,Z_h$ and the numbers $N_h,\delta_h,M_h$ are independent of
$z,\theta,t$. By \eqref{2026093003},
$|t\eps\phi(p-q)|\le a\langle p-q\rangle^{-r}$ for
$p,q\in\mathbb Z^d$ and $t\in[0,1]$.

\subsection{A common choice of constants}\label{sec:induction-proof}
The estimates below require three smallness conditions, one for
elimination, one for eigenvalue drift and one for eigenvector
reconstruction. We fix all coefficients before stating the induction
assertions.

Let $\mathrm C_{5}=\mathrm C_{5}(d,r,\tau)$ be a finite constant such that
\begin{equation}\label{2026100904}
 \mathrm C_{5}\ge32\cdot16^{2r}\mathrm C_{3}
\end{equation}
and the diagonal comparison \eqref{2026101304} holds with coefficient
$\mathrm C_{5}/2$. Under the stated matrix hypotheses,
Proposition~\ref{prop:schur-comparison} and
Corollary~\ref{cor:diagonal-comparison} prove \eqref{2026101304}
with a finite coefficient depending only on $d,r,\tau$,
so $\mathrm C_{5}$ can be fixed before the induction.
The factor in \eqref{2026100904} ensures that the
last Schur elimination in \eqref{2026101601}
also has coefficient at most $\mathrm C_{5}/2$.

For the fixed $0<\rho<r-\tau$, put
\begin{equation}\label{2026100905}
 G(s)=1+\frac{16^\tau}{1-16^{\tau+s-r}},\qquad 0\le s<r-\tau.
\end{equation}
The scales in \eqref{2026100605}--\eqref{2026100606} give
\[
 \sum_{h\ge0}\delta_h^{-1}M_h^{s-r}=\beta^{-1}G(s),
 \qquad 0\le s<r-\tau.
\]
Fix a finite constant
$\mathrm C_{6}=\mathrm C_{6}(d,r,\tau,\gamma,L,\rho)$ satisfying
\begin{equation}\label{2026100906}
 \mathrm C_{6}\ge2\mathrm C_{4}\beta^{-1}G(\rho),\qquad
 \mathrm C_{6}\ge\frac{2^{\tau+4}\mathrm C_{3}}{L\gamma}.
\end{equation}
The first inequality controls the product of the reconstruction
matrices; the second controls the use of eigenvalue separation in
\eqref{2026101801}.

We require
\begin{align}
 &\mathrm C_{5} a\sum_{h\ge0}\delta_h^{-1}M_h^{-r}\le\frac1{64},
 \label{2026100907}\\
 &\frac{\mathrm C_{5} a^2\beta^{-2}}{1-16^{\tau-2r}}\le\frac1{16},
 \label{2026100908}\\
 &\mathrm C_{6} a\le\frac14.
 \label{2026100909}
\end{align}
Choose a threshold
\begin{equation}\label{2026101001}
 \eps_0=\eps_0(d,r,\tau,\gamma,L,\rho,\mathrm C_{0})>0
\end{equation}
so that \eqref{2026100907}--\eqref{2026100909}
hold with $a=\mathrm C_{0}\eps_0$.
The left-hand sides tend to zero as $a\to0$, so such a threshold exists.
Since $a\le\mathrm C_{0}|\eps|$, the same three inequalities hold
whenever $|\eps|<\eps_0$.
The dependence in \eqref{2026101001} agrees with
Theorem~\ref{thm:main}, since $\mathrm C_{0}=\mathrm C_{0}(r,\phi)$.
Moreover, \eqref{2026100907} implies
\[
 \sum_{h\ge0}\eta_h\le(64\mathrm C_{5})^{-1}
 \le(64\mathrm C_{3})^{-1},\qquad
 a/\beta\le(2\mathrm C_{3})^{-1}.
\]
Thus the comparison estimates and the initial construction apply.
The constants $\mathrm C_{3},\ldots,\mathrm C_{6}$ and $\eps_0$ remain
fixed throughout the induction and the localization proof,
independently of $n$, $\theta$ and $t$.
Finally, set
\begin{equation}\label{2026101002}
 \kappa=\frac{\rho^2}{\tau+\rho}>0.
\end{equation}

\subsection{The induction assertions}
We state assertions \textup{(1)--(5)} at $t=1$ and require
validity for every $t\in[0,1]$ in assertion \textup{(6)}.

\begin{theorem}[Multiscale estimates]\label{thm:finite-volume-induction}
Assume \eqref{2026092601}, the potential assumptions in Section~\ref{sec:model}
and $\om\in\dc$, with $r>\tau>d$.
Under \eqref{2026100907},
\eqref{2026100908} and \eqref{2026100909},
there are finite sets $(B_m)_{m\ge0}$, with $B_0=\{o\}$,
for which the functions
\begin{equation}
 E_0(\theta;t)=f(\theta),\qquad E_m(\theta;t)=\lambda(B_m;\theta;t),
 \qquad m\ge1,\quad\theta\in\Theta,\quad t\in[0,1],
 \label{2026101003}
\end{equation}
satisfy the following assertions at every step $n\ge1$.
Here $\lambda(B_m;\theta;t)$ is the branch defined in
\eqref{2026100403} for
$\cH_{B_m}(\theta;t)$, and we write $E_m=E_m(\cdot;1)$.
In particular, each $E_m(\cdot;t)$ is real and $1$-periodic.
\begin{enumerate}[label=\textup{(\arabic*)},ref=\arabic*,leftmargin=*]
\item\label{item:induction-geometry}
\textup{\textbf{Geometry.}}
For $1\le m\le n$,
\begin{align}
 B_{m-1}&\subset B_m,\qquad
 \Lambda_{4N_m}\subset B_m\subset\Lambda_{8N_m},
 \label{2026101004}\\
 p\in B_m\cap Z_h,\quad0\le h\le m-2
 &\quad\Longrightarrow\quad B_{h+1}(p)+\Lambda_{N_{h+1}}\subset B_m,
 \label{2026101005}\\
 B_m\cap Z_{m-1}&=\{o\}.
 \label{2026101006}
\end{align}
The sets $\Lambda_R$ and the addition of sets are defined in
\eqref{2026092705} and \eqref{2026092703}.
For $0\le h\le n-1$, we also have
\begin{equation}\label{2026101007}
 p\in Z_h^\circ
 \quad\Longrightarrow\quad o\notin B_{h+1}(p)+\Lambda_{N_{h+1}}.
\end{equation}

\item\label{item:induction-spectrum}
\textup{\textbf{Local spectrum.}}
For $1\le m\le n$, $\theta\in\Theta$ and $z\in J_m(\theta)$,
$\cA_{B_m^\circ}(z;\theta)$ is invertible.
The function $s_m(\cdot;\theta)$ has exactly one zero in $J_m(\theta)$,
namely $E_m(\theta)$, and satisfies
\begin{equation}\label{2026101008}
 |s_m(z;\theta)|\ge|z-E_m(\theta)|,
 \qquad z\in J_m(\theta).
\end{equation}
The eigenvalue $E_m(\theta)$ of $\cH_{B_m}(\theta)$ is simple and
\begin{equation}\label{2026101009}
 \dist\!\left(E_m(\theta),
       \Spec\cH_{B_m}(\theta)\setminus\{E_m(\theta)\}\right)
 \ge3\delta_{m-1}.
\end{equation}

\item\label{item:induction-branches}
\textup{\textbf{Eigenvalue estimates.}}
With the coefficient $\mathrm C_{5}$ fixed in Subsection~\ref{sec:induction-proof},
\begin{align}
 &\|E_m-E_{m-1}\|_\infty
 \le\mathrm C_{5} a^2\beta^{-1}N_{m-1}^{\tau-2r},\quad 1\le m\le n,
 \label{2026101101}\\
 &\|E_m-E_\ell\|_\infty\le\delta_\ell/8,\quad 0\le\ell<m\le n.
 \label{2026101102}
\end{align}
Here $\|\cdot\|_\infty$ is the supremum norm specified in
Theorem~\ref{thm:main}.
For $1\le m\le n$, we also have
\begin{equation}\label{2026101103}
 E_m(\theta_2)-E_m(\theta_1)\ge L(\theta_2-\theta_1),
 \qquad \theta_1<\theta_2,\quad \theta_1,\theta_2\in\Theta\cap[0,1).
\end{equation}

\item\label{item:induction-matrices}
\textup{\textbf{Elimination estimates.}}
For $0\le i\le n-1$, $\theta\in\Theta$, $z\in\mathbb R$ and
$\Omega\in\Upsilon_i$, the matrices
$(\cA_\Omega^{(h-1)})_{Q_h\cap\Omega}$ are invertible for $0\le h\le i$.
Moreover,
\begin{align}
 &|\cA_\Omega^{(h-1)}(p,p)|\ge\delta_h/2,\quad p\in Q_h\cap\Omega,\quad0\le h\le i,
 \label{2026101104}\\
 &|\cA_\Omega^{(h-1)}(p,q)|\le2a\langle p-q\rangle^{-r},\quad p,q\in S_{h-1}\cap\Omega,\quad p\ne q,\quad 0\le h\le i+1.
 \label{2026101105}
\end{align}
Estimate \eqref{2026101105} includes
$\cA_\Omega^{(i)}$, the matrix after eliminating $I_i$.

\item\label{item:induction-vectors}
\textup{\textbf{Eigenvectors and increments.}}
For $0\le m\le n$ and $\theta\in\Theta$, there is a vector
$v_{m,\theta}$ supported on $B_m$, with $v_{0,\theta}=e_o$, such that
\begin{equation}\label{2026101106}
 \cH_{B_m}(\theta)\left.v_{m,\theta}\right|_{B_m}
 =E_m(\theta)\left.v_{m,\theta}\right|_{B_m},\qquad
 v_{m,\theta}(o)=1.
\end{equation}
Extend $v_{m,\theta}$ by zero to $\mathbb Z^d$. Then
\begin{align}
 &\|v_{m,\theta}-e_o\|_{1,\rho}
 \le\mathrm C_{6} a\le1/4,\quad 0\le m\le n,
 \label{2026101201}\\
 &\sup_{\theta\in\Theta}
 \|v_{m,\theta}-v_{m-1,\theta}\|_1
 \le12\mathrm C_{6}aN_{m-1}^{-\kappa},\quad 1\le m\le n.
 \label{2026101202}
\end{align}
The norm $\|\cdot\|_{1,\rho}$ is defined in
\eqref{2026100306}; $\mathrm C_{6}$ and $\kappa$ are fixed in
\eqref{2026100906} and \eqref{2026101002}.

\item\label{item:induction-homotopy}
\textup{\textbf{Uniformity for $t\in[0,1]$.}}
For every $t\in[0,1]$, assertions \textup{(1)--(5)} hold for
\[
 \cH(\cdot;t),\quad E_m(\cdot;t),\quad s_m(\cdot;\cdot;t),
 \quad J_m(\cdot;t),\quad v_{m,\theta}(t),
\]
with $B_m$, $N_h$, $\delta_h$, $M_h$ and all constants independent of $t$.
In assertion \textup{(4)}, use $S_h(z;\theta;t)$,
$Q_h(z;\theta;t)$, $\Upsilon_i(z;\theta;t)$ and
$\cA_\Omega^{(h)}(z;\theta;t)$ as defined after
\eqref{2026100903}.
For $0\le m\le n$ and $\theta\in\Theta$, equations
\eqref{2026101003} and \eqref{2026100404} give
\begin{equation}\label{2026101203}
 E_m(\theta;0)=f(\theta),\qquad
 E_m(\theta;1)=E_m(\theta).
\end{equation}
By \eqref{2026092905}, $t\mapsto E_m(\theta;t)$ is continuous.
Moreover, $v_{m,\theta}(0)=e_o$ and
$t\mapsto v_{m,\theta}(t)$ is continuous in $\ell^2$.
\end{enumerate}
\end{theorem}

We record how assertion \textup{(1)} gives the regions required in
assertion \textup{(4)}. Fix $1\le m\le n$, $\theta\in\Theta$,
$z\in\mathbb R$ and $t\in[0,1]$. If
$S_h(z;\theta;t)\subset Z_h$ for $0\le h<m$, then
\eqref{2026101005} gives
$B_m\in\Upsilon_{m-1}(z;\theta;t)$.
For $p\in B_m^\circ\cap S_h$ with $0\le h\le m-2$,
\eqref{2026101007} also gives $o\notin B_{h+1}(p)$;
hence $B_m^\circ\in\Upsilon_{m-1}$.
Moreover, \eqref{2026101005} gives
$\dist(x,B_m^c)\ge N_{h+1}$ for $p\in B_m\cap S_h$,
$0\le h\le m-2$ and $x\in B_{h+1}(p)$.
Finally, by \eqref{2026101006}, the remaining index set after
eliminating $Q_0,\ldots,Q_{m-1}$ is
$B_m\cap S_{m-1}\subset\{o\}$.

We verify Theorem~\ref{thm:finite-volume-induction} at $n=1$.
By \eqref{2026100105}, $B_1=\Lambda_{4N_1}$ and $B_0=\{o\}$.
Equation \eqref{2026100703} gives $|p|>33600$ for $p\in Z_0^\circ$.
Since $4N_1=64$ and $5N_1=80$, we have
\[
 B_1\cap Z_0=\{o\},\qquad
 o\notin B_1(p)+\Lambda_{N_1}\quad(p\in Z_0^\circ).
\]
The range $0\le h\le m-2$ in \eqref{2026101005} is empty for $m=1$.
This proves assertion \textup{(1)}.
Proposition~\ref{prop:initial-step} and \eqref{2026100504}
prove assertion \textup{(2)}.
By \eqref{2026100303} and
$\mathrm C_{5}\ge2\mathrm C_{3}$, the bounds
\eqref{2026101101} and \eqref{2026101102} hold for $m=1$.
Equation \eqref{2026100502} and
Proposition~\ref{prop:kps-label} give \eqref{2026101003}
and \eqref{2026101103}, completing assertion \textup{(3)}.
Since $\Upsilon_0$ consists of all finite sets,
Lemma~\ref{lem:initial-elimination}, \eqref{2026100104}
and \eqref{2026100103} prove assertion \textup{(4)} at $i=0$.

The eigenvector in \eqref{2026100305} satisfies
\eqref{2026101106}.
By \eqref{2026100401} and
\eqref{2026100906},
\[
 \|v_{1,\theta}(t)-e_o\|_{1,\rho}
 \le\frac{\mathrm C_{4}a}{4\beta}\le\mathrm C_{6} a\le1/4.
\]
Since $N_0=1$, the same bound proves
\eqref{2026101202} for $m=1$.
This proves assertion \textup{(5)}.
The first construction holds for every $t\in[0,1]$, and
\eqref{2026100302} and \eqref{2026100305}
give the continuity and endpoint identities in assertion \textup{(6)}.

\section{The induction step}\label{sec:induction-step}
Assume that the assertions of Theorem~\ref{thm:finite-volume-induction}
hold at step $n\ge1$.
We follow the equations in Subsection~\ref{sec:next-equations}:
compare $\cA_\Omega^{(n-1)}(p,p)$ with
$s_n(z;\theta+p\cdot\om)$, solve for
the coordinates in $Q_n$, and construct the next local eigenvalue
and corresponding eigenvector on a fixed region $B_{n+1}$.
We also estimate $E_{n+1}-E_n$ and $v_{n+1,\theta}-v_{n,\theta}$.
The case $n=1$ is included in the argument.
Fix $t\in[0,1]$ and omit $t$ from the notation except when checking the homotopy.

\subsection{Comparison on two regions}\label{sec:entrywise-comparison}
Fix $\theta\in\Theta$, $z\in\mathbb R$ and a finite
$F\in\Upsilon_n$, and consider
\[
 \cA_F^{(n-1)}\left.v\right|_{S_{n-1}\cap F}=0.
\]
To solve for $\left.v\right|_{Q_n\cap F}$, we compare
$\cA_F^{(n-1)}(p,p)$ with $s_n(z;\theta+p\cdot\om)$ for
$p\in Q_n\cap F$.
By \eqref{2026100804}, $B_n(p)\subset F$.
Thus the application uses $(\Omega,\Omega')=(B_n(p),F)$
in the comparison below.
The arguments $(z,\theta,t)$ are fixed while the finite region varies.

We first compare the matrices after one elimination.

\begin{lemma}\label{lem:first-region-comparison}
Let $\Omega\subset\Omega'$ be finite, let $Q\subset\mathbb Z^d$,
and let $\cM$ be a Hermitian matrix on $\Omega'$ such that
\[
 |\cM(x,y)|\le2a\langle x-y\rangle^{-r}\quad(x\ne y),
 \qquad |\cM(x,x)|\ge\beta/2\quad(x\in Q\cap\Omega').
\]
Assume $a/\beta\le(64\mathrm C_{3})^{-1}$ and put
$P=\Omega\setminus Q$, $Y=Q\cap(\Omega'\setminus\Omega)$.
The Schur complements
\[
 \mathcal B_\Omega=\Sch_{Q\cap\Omega}(\cM_\Omega),\qquad
 \mathcal B_{\Omega'}=\Sch_{Q\cap\Omega'}(\cM)
\]
exist and, for $p,q\in P$, satisfy
\begin{equation}\label{2026101204}
 |\mathcal B_{\Omega'}(p,q)-\mathcal B_\Omega(p,q)|
 \le72a^2\beta^{-1}\sum_{u\in Y}
              \langle p-u\rangle^{-r}\langle q-u\rangle^{-r}.
\end{equation}
\end{lemma}
\begin{proof}
	By the assumption $a/\beta\le(64\mathrm C_{3})^{-1}$,
	$4a/\beta\le(16\mathrm C_{3})^{-1}<(2\mathrm C_{3})^{-1}$.
	Since $\Omega'$ is $1$-separated, the assumptions on $\cM$
	verify the hypotheses of Lemma~\ref{lem:sparse-schur} with
	$\fX=\Omega'$, $\fP=\Omega'\setminus(Q\cap\Omega)$,
	$\fQ=Q\cap\Omega$, $M=1$, $a_*=2a$ and $\Delta=\beta/2$.
	Hence $\Sch_{Q\cap\Omega}(\cM)$ exists.
	By \eqref{2026092803},
	\[
	\bigl(\Sch_{Q\cap\Omega}(\cM)\bigr)_P
	=\Sch_{Q\cap\Omega}(\cM_\Omega)
	=\mathcal B_\Omega.
	\]
	
	The assumption $a/\beta\le(64\mathrm C_{3})^{-1}$
	implies
	\[
	16\mathrm C_{3}a/\beta\le1,
	\qquad
	16\mathrm C_{3}(a/\beta)^2\le\frac14.
	\]
	Thus \eqref{2026102101} gives, for distinct
	$x,y\in\Omega'\setminus(Q\cap\Omega)$,
	\[
	|\Sch_{Q\cap\Omega}(\cM)(x,y)|
	\le
	\bigl(2+16\mathrm C_{3}a/\beta\bigr)
	a\langle x-y\rangle^{-r}
	\le3a\langle x-y\rangle^{-r}.
	\]
	For $u\in Y$, \eqref{2026102102} gives
	\[
	|\Sch_{Q\cap\Omega}(\cM)(u,u)|
	\ge\beta/2-16\mathrm C_{3}a^2/\beta
	\ge\beta/4.
	\]
	
	The bounds on the entries of $\Sch_{Q\cap\Omega}(\cM)$ and
	$12a/\beta\le(2\mathrm C_{3})^{-1}$ verify the assumptions of
	Lemma~\ref{lem:sparse-schur} with
	$\fX=\Omega'\setminus(Q\cap\Omega)$, $\fP=\Omega'\setminus Q$,
	$\fQ=Y$, $M=1$, $a_*=3a$ and $\Delta=\beta/4$.
	For $u\in Y$ and $q\in P$, \eqref{2026102103} yields
	\[
	\left|
	\left[
	\bigl((\Sch_{Q\cap\Omega}(\cM))_Y\bigr)^{-1}
	(\Sch_{Q\cap\Omega}(\cM))_{Y,q}
	\right](u)
	\right|
	\le24a\beta^{-1}\langle u-q\rangle^{-r}.
	\]
	Because $(Q\cap\Omega)\sqcup Y=Q\cap\Omega'$,
	\eqref{2026092805} implies
	\[
	\mathcal B_{\Omega'}
	=\Sch_Y\bigl(\Sch_{Q\cap\Omega}(\cM)\bigr).
	\]
	Hence, for $p,q\in P$, \eqref{2026092803} gives
	\begin{equation*}
 |\mathcal B_{\Omega'}(p,q)-\mathcal B_\Omega(p,q)|
		\le
		\sum_{u\in Y}
		3a\langle p-u\rangle^{-r}
		\,24a\beta^{-1}\langle u-q\rangle^{-r} =72a^2\beta^{-1}
		\sum_{u\in Y}
		\langle p-u\rangle^{-r}\langle q-u\rangle^{-r}.
\end{equation*}
	This proves \eqref{2026101204}.
	\end{proof}

After further eliminations, the matrices on the common indices may
already differ. The next proposition estimates the accumulated difference.

\begin{proposition}\label{prop:schur-comparison}
Let $\Omega\subset\Omega'$ be finite and let $i\ge1$.
For this proposition, let
$S_{-1}=\mathbb Z^d\supset S_0\supset\cdots\supset S_{i-1}$
be any decreasing family of sets.
Define $Q_k$ by \eqref{2026100706} for $0\le k<i$.
Assume that
\[
 |x-y|\ge M_k,\qquad x,y\in S_{k-1},\quad x\ne y,\quad 0\le k<i.
\]
For $F\in\{\Omega,\Omega'\}$, define $\cA_F^{(k)}$ by
\eqref{2026100708} using $S_{-1},\ldots,S_{i-1}$.
Suppose that, for each $0\le k<i$, the matrices
$\cA_{\Omega}^{(k-1)}$ and $\cA_{\Omega'}^{(k-1)}$ exist and satisfy
the following conditions:
\begin{enumerate}[label=\textup{(\arabic*)}]
\item For $F\in\{\Omega,\Omega'\}$,
      \[
       |\cA_{F}^{(k-1)}(x,y)|\le2a\langle x-y\rangle^{-r},
       \qquad x,y\in S_{k-1}\cap F,\quad x\ne y.
      \]
\item For $F\in\{\Omega,\Omega'\}$,
      \[
       |\cA_{F}^{(k-1)}(x,x)|\ge\delta_k/2,
       \qquad x\in Q_k\cap F.
      \]
\item
      \[
       \dist(x,\Omega^c)\ge M_k,
       \qquad x\in S_{k-1}\cap\Omega.
      \]
\end{enumerate}
If
\begin{equation}\label{2026101205}
 a\sum_{k\ge0}\delta_k^{-1}M_k^{-r}\le(64\mathrm C_{3})^{-1},
\end{equation}
then, for $p,q\in S_{i-1}\cap\Omega$,
\begin{equation}
 |\cA_{\Omega'}^{(i-1)}(p,q)-\cA_{\Omega}^{(i-1)}(p,q)| \le32\mathrm C_{3}
 a^2\sum_{k=0}^{i-1}\delta_k^{-1}
 \sum_{u\in Q_k\cap(\Omega'\setminus\Omega)}
 \langle p-u\rangle^{-r}\langle q-u\rangle^{-r}.
 \label{2026101206}
\end{equation}
\end{proposition}
\begin{proof}
We proceed by induction on the number of eliminated sets.
Write
\[
 Y_k=Q_k\cap(\Omega'\setminus\Omega),\qquad 0\le k<i.
\]
Use $I_k$ and $P_k$ from \eqref{2026100707}
for the region $\Omega$.
By \eqref{2026100708},
\[
 (\cA_{\Omega'}^{(-1)})_\Omega=\cA_\Omega^{(-1)}.
\]
Since $\delta_0=\beta$, conditions \textup{(1)--(2)} and
\eqref{2026101205} verify the assumptions of
Lemma~\ref{lem:first-region-comparison} with $Q=Q_0$ and
$\cM=\cA_{\Omega'}^{(-1)}$.
Thus \eqref{2026101204} holds.

\emph{The general step.}
Let $1\le k<i$ and assume that, for $p,q\in P_{k-1}$,
\begin{align}
 |\cA_{\Omega'}^{(k-1)}(p,q)-\cA_{\Omega}^{(k-1)}(p,q)|
 &\le16\mathrm C_{3}a^2\sum_{h=0}^{k-1}\delta_h^{-1}
       \prod_{h<\ell<k}(1+16\mathrm C_{3}\eta_\ell)^2\notag\\*
 &\qquad\times\sum_{u\in Y_h}
                   \langle p-u\rangle^{-r}\langle q-u\rangle^{-r}.
 \label{2026101301}
\end{align}
An empty product equals $1$.
Since $72\le16\mathrm C_{3}$, \eqref{2026101204} proves
\eqref{2026101301} for $k=1$.
At step $k$, we first eliminate $I_k$ from the common indices
$P_{k-1}=P_k\sqcup I_k$ and then $Y_k$ in $\Omega'$.
Write $\cM=\cA_{\Omega}^{(k-1)}$ and
$\widetilde{\cM}=\cA_{\Omega'}^{(k-1)}$.
Unlike the original matrices at step $0$,
$\widetilde{\cM}_{P_{k-1}}$ and $\cM$ need not agree.
We first estimate
$\Sch_{I_k}(\widetilde\cM_{P_{k-1}})-\Sch_{I_k}(\cM)$,
then eliminate $Y_k$ in $\Omega'$.

\emph{Estimating $\Sch_{I_k}(\widetilde\cM_{P_{k-1}})-\Sch_{I_k}(\cM)$.}
For $F\in\{\Omega,\Omega'\}$, conditions \textup{(1)--(2)} and
\eqref{2026101205} verify the assumptions of
Lemma~\ref{lem:sparse-schur} for $\cA_F^{(k-1)}$ with
$\fX=S_{k-1}\cap F$, $\fP=(S_{k-1}\cap F)\setminus I_k$,
$\fQ=I_k$, $M=M_k$, $a_*=2a$ and $\Delta=\delta_k/2$.
Indeed, $S_{k-1}$ is $M_k$-separated and
$4\eta_k\le4(64\mathrm C_{3})^{-1}\le(2\mathrm C_{3})^{-1}$.
Thus $\cM_{I_k}$ and $\widetilde\cM_{I_k}$ are invertible.
Let $\cW=\cW(\cM)$ and
$\widetilde\cW=\cW(\widetilde\cM_{P_{k-1}})$, where the reconstruction
matrices in \eqref{2026092804} use $P=P_k$ and $Q=I_k$.
By \eqref{2026092806}, we have
\[
 \Sch_{I_k}(\widetilde\cM_{P_{k-1}})-\Sch_{I_k}(\cM)
 =\cW^*(\widetilde{\cM}_{P_{k-1}}-\cM)\widetilde{\cW}.
\]

Fix $u\in Y_h$ with $0\le h<k$ and $p\in P_k$.
Since $u\notin\Omega$, condition \textup{(3)} implies
\[
 \dist(u,P_{k-1})\ge M_k.
\]
The set $P_{k-1}\subset S_{k-1}$ is $M_k$-separated and
$I_k\subset P_{k-1}\setminus\{p\}$.
Thus \eqref{2026102103} and \eqref{2026101906} imply
\begin{align*}
 &\sum_{x\in P_{k-1}}|\cW(x,p)|\langle x-u\rangle^{-r}
 \le \langle p-u\rangle^{-r}
       +8a\delta_k^{-1}\sum_{x\in I_k}\langle p-x\rangle^{-r}\langle x-u\rangle^{-r}\\
 &\le(1+16\mathrm C_{3}\eta_k)\langle p-u\rangle^{-r}.
\end{align*}
By \eqref{2026102103} and \eqref{2026101906}, we also have
\[
 \sum_{y\in P_{k-1}}|\widetilde\cW(y,q)|\langle y-u\rangle^{-r}
 \le(1+16\mathrm C_{3}\eta_k)\langle q-u\rangle^{-r},\qquad q\in P_k.
\]
Multiplying the bounds for $\cW$ and $\widetilde\cW$, we obtain
\begin{equation*}
 \sum_{x,y\in P_{k-1}}|\cW(x,p)|
       \langle x-u\rangle^{-r}\langle y-u\rangle^{-r}|\widetilde{\cW}(y,q)| \le(1+16\mathrm C_{3}\eta_k)^2\langle p-u\rangle^{-r}\langle q-u\rangle^{-r}.
\end{equation*}

\emph{Eliminating $Y_k$.}
Let
$\widehat{\cM}^{(k)}=\Sch_{I_k}(\widetilde{\cM})$.
By \eqref{2026102102},
\[
 |\widehat{\cM}^{(k)}(u,u)-\widetilde{\cM}(u,u)|
 \le16\mathrm C_{3}a^2\delta_k^{-1}M_k^{-2r}
 =16\mathrm C_{3}\eta_k^2\delta_k,\qquad u\in Y_k.
\]
By Proposition~\ref{prop:schur-comparison}\,\textup{(2)},
\eqref{2026102101} and $\eta_k\le(64\mathrm C_{3})^{-1}$, we obtain
\begin{align*}
 |\widehat\cM^{(k)}(u,u)|&\ge\delta_k/2-16\mathrm C_{3}\eta_k^2\delta_k
 \ge\delta_k/4,\qquad u\in Y_k,\\
 |\widehat\cM^{(k)}(x,y)|&\le(2+16\mathrm C_{3}\eta_k)a\langle x-y\rangle^{-r}
 \le3a\langle x-y\rangle^{-r},\qquad x\ne y.
\end{align*}
In the second line, $x,y\in(S_{k-1}\cap\Omega')\setminus I_k$.
The entry bounds on $\widehat{\cM}^{(k)}$ and
$12\eta_k\le12(64\mathrm C_{3})^{-1}\le(2\mathrm C_{3})^{-1}$
verify the assumptions of Lemma~\ref{lem:sparse-schur} with
$\fX=(S_{k-1}\cap\Omega')\setminus I_k$,
$\fP=S_k\cap\Omega'$, $\fQ=Y_k$, $M=M_k$,
$a_*=3a$ and $\Delta=\delta_k/4$.
Hence $\widehat{\cM}^{(k)}_{Y_k}$ is invertible.
By \eqref{2026102103}, we have
\[
 |((\widehat{\cM}^{(k)}_{Y_k})^{-1}\widehat{\cM}^{(k)}_{Y_k,q})(u)|
 \le24a\delta_k^{-1}\langle u-q\rangle^{-r},
 \quad u\in Y_k,\quad q\in P_k.
\]
It follows that
\begin{equation}\label{2026101302}
 |\widehat{\cM}^{(k)}_{p,Y_k}(\widehat{\cM}^{(k)}_{Y_k})^{-1}\widehat{\cM}^{(k)}_{Y_k,q}|
 \le72a^2\delta_k^{-1}
       \sum_{u\in Y_k}\langle p-u\rangle^{-r}\langle q-u\rangle^{-r}.
\end{equation}

\emph{Summation over the steps.}
Since $I_k\sqcup Y_k=Q_k\cap\Omega'$, we have
\[
 \Sch_{Y_k}\bigl(\Sch_{I_k}(\widetilde\cM)\bigr)
 =\Sch_{Q_k\cap\Omega'}(\widetilde\cM)
 =\cA_{\Omega'}^{(k)}.
\]
Substituting \eqref{2026101301} into
\eqref{2026092806} and adding
\eqref{2026101302} proves
\eqref{2026101301} with $k$ replaced by $k+1$.
Moreover, \eqref{2026101205} implies
\begin{equation}\label{2026101303}
 \prod_{k<i}(1+16\mathrm C_{3}\eta_k)^2
 \le\exp\!\left(32\mathrm C_{3}\sum_{k\ge0}\eta_k\right)
 \le e^{1/2}<2.
\end{equation}
Since $Y_h=Q_h\cap(\Omega'\setminus\Omega)$, combining
\eqref{2026101303} with
\eqref{2026101301} at $k=i$ proves
\eqref{2026101206}.
\end{proof}

\begin{corollary}[Diagonal comparison]\label{cor:diagonal-comparison}
Under the assumptions of Proposition~\ref{prop:schur-comparison},
let $p\in S_{i-1}\cap\Omega$ and let
$R=\dist(p,\Omega^c)\ge1$.
If $M_k\le R$ for $0\le k<i$, then
\begin{equation}\label{2026101304}
 |\cA_{\Omega'}^{(i-1)}(p,p)-\cA_{\Omega}^{(i-1)}(p,p)|
 \le\frac{\mathrm C_{5}}{2}a^2\beta^{-1}R^{\tau-2r}.
\end{equation}

\end{corollary}
\begin{proof}
For $0\le k<i$, $Q_k\subset S_{k-1}$ is $M_k$-separated and
\[
 |p-u|\ge\dist(p,\Omega^c)=R\ge M_k,
 \qquad u\in Q_k\cap(\Omega'\setminus\Omega).
\]
By \eqref{2026101907},
\[
 \sum_{u\in Q_k\cap(\Omega'\setminus\Omega)}
 \langle p-u\rangle^{-2r}\le\mathrm C_{3}M_k^{-d}R^{d-2r}.
\]
By \eqref{2026101206},
\[
 |\cA_{\Omega'}^{(i-1)}(p,p)-\cA_{\Omega}^{(i-1)}(p,p)|
 \le32\mathrm C_{3}^2
 a^2R^{d-2r}\sum_{k<i}\delta_k^{-1}M_k^{-d}.
\]
Since $\tau>d$ and $M_k\le R$, we have
\begin{equation*}
 \sum_{k<i}\delta_k^{-1}M_k^{-d}
 =\beta^{-1}\left(1+16^\tau\sum_{k=1}^{i-1}N_{k-1}^{\tau-d}\right) \le\beta^{-1}\left(1+
       \frac{16^\tau R^{\tau-d}}{1-16^{d-\tau}}\right)
 \le\left(1+\frac{16^\tau}{1-16^{d-\tau}}\right)\beta^{-1}R^{\tau-d}.
\end{equation*}
The sum over $k\ge1$ is zero when $i=1$.
The coefficient
$32\mathrm C_{3}^2(1+16^\tau/(1-16^{d-\tau}))$ is finite and depends
only on $d,r,\tau$. This is the coefficient required to be at most
$\mathrm C_{5}/2$ when $\mathrm C_{5}$ was fixed in
Subsection~\ref{sec:induction-proof}, and proves \eqref{2026101304}.
\end{proof}

Dividing \eqref{2026101304} by $\beta R^{-\tau}$, we obtain
\begin{equation}\label{2026101305}
 \frac{|\cA_{\Omega'}^{(i-1)}(p,p)-\cA_{\Omega}^{(i-1)}(p,p)|}
      {\beta R^{-\tau}}
 \le\frac{\mathrm C_{5}}{2}a^2\beta^{-2}R^{2\tau-2r}.
\end{equation}
Since $r>\tau$, the right-hand side of
\eqref{2026101305} is at most
$\frac{\mathrm C_{5}}{2}a^2\beta^{-2}$ for every $R\ge1$.
Thus \eqref{2026100908} bounds the quotient in
\eqref{2026101305} uniformly for $R\ge1$.

\subsection{The new diagonal bound and the next elimination}
\label{sec:general-inductive-step}
We use \eqref{2026101103} only for $1\le h\le n$;
for $h=0$, the required inequality follows from $E_0=f$ and
\eqref{2026092604}.

\begin{lemma}\label{lem:next-elimination}
For every $\theta\in\Theta$, $z\in\mathbb R$ and
$\Omega\in\Upsilon_n$, the matrix
$(\cA_\Omega^{(n-1)})_{Q_n\cap\Omega}$ is invertible and
\begin{align}
 &|\cA_\Omega^{(n-1)}(p,p)|\ge\delta_n/2,\quad p\in Q_n\cap\Omega,
 \label{2026101306}\\
 &|\cA_\Omega^{(n)}(p,q)|\le2a\langle p-q\rangle^{-r},\quad p,q\in S_n\cap\Omega,\quad p\ne q.
 \label{2026101307}
\end{align}
Consequently, assertion \textup{(4)} extends to $i=n$.
\end{lemma}
\begin{proof}
Fix $\theta\in\Theta$, $z\in\mathbb R$, $\Omega\in\Upsilon_n$
and $p\in S_{n-1}\cap\Omega$.
By \eqref{2026100804}, $B_n(p)\subset\Omega$.
If $1\le h\le n$ and $q\in B_n(p)\cap S_{h-1}$, then
$p,q\in\Res_{h-1}$ by \eqref{2026100705}.
By \eqref{2026100507},
\begin{equation}\label{2026101401}
 L\|(q-p)\cdot\om\|_{\mathbb T}
 \le|E_{h-1}(\theta+q\cdot\om)-E_{h-1}(\theta+p\cdot\om)|
 <2\delta_{h-1}.
\end{equation}
By \eqref{2026101401} and \eqref{2026100702},
$q-p\in Z_{h-1}$.
For $1\le h\le n-1$, \eqref{2026101005} implies
$B_h(q)\subset B_n(p)$.  This proves
$B_n(p)\in\Upsilon_{n-1}$.
For $h=n$, \eqref{2026101006} implies $q=p$.
Since $p\in S_{n-1}$, we have
\begin{equation}\label{2026101402}
 B_n(p)\cap S_{n-1}=\{p\}.
\end{equation}

We verify the assumptions of Proposition~\ref{prop:schur-comparison}
with $i=n$ and regions $B_n(p)\subset\Omega$.
Since $B_n(p),\Omega\in\Upsilon_{n-1}$,
\eqref{2026101104}--\eqref{2026101105}
verify conditions \textup{(1)--(2)} through step $n-1$.
For $1\le k<n$ and $x\in B_n(p)\cap S_{k-1}$,
\eqref{2026100804} implies
\[
 B_k(x)\subset B_n(p),\qquad
 \dist(x,B_n(p)^c)\ge4N_k\ge M_k.
\]
At $k=0$, $\dist(x,B_n(p)^c)\ge1=M_0$ for every
$x\in B_n(p)$.  This proves condition \textup{(3)}.
The sets $S_{k-1}\subset\Res_{k-1}$ are $M_k$-separated for $1\le k<n$:
indeed, \eqref{2026100507} and \eqref{2026100704} give
\[
 |x-y|\ge\frac12\left(\frac{L\gamma}{2\beta}\right)^{1/\tau}N_{k-1}
 >M_k\qquad(x,y\in S_{k-1},\ x\ne y).
\]
The set $S_{-1}=\mathbb Z^d$ is $M_0$-separated.
Finally, \eqref{2026100907} verifies
\eqref{2026101205}.
By \eqref{2026101402}, eliminating
$Q_0\cap B_n(p),\ldots,Q_{n-1}\cap B_n(p)$ leaves the index set $\{p\}$.
Thus \eqref{2026092805} and \eqref{2026092708} imply
\[
 \cA_{B_n(p)}^{(n-1)}(p,p)=s_n(z;\theta+p\cdot\om).
\]
Here $z\in J_n(\theta+p\cdot\om)$ because $p\in S_{n-1}$ implies
$|z-E_{n-1}(\theta+p\cdot\om)|<\delta_{n-1}$.
Since $\dist(p,B_n(p)^c)\ge4N_n\ge M_k$ for $0\le k<n$,
Corollary~\ref{cor:diagonal-comparison} gives
\begin{equation}\label{2026101403}
 |\cA_\Omega^{(n-1)}(p,p)-s_n(z;\theta+p\cdot\om)|
 \le\frac{\mathrm C_{5}}{2}a^2\beta^{-1}N_n^{\tau-2r}.
\end{equation}
For $p\in Q_n=S_{n-1}\setminus\Res_n$,
\eqref{2026101008} gives
\begin{equation}\label{2026101404}
 |s_n(z;\theta+p\cdot\om)|
 \ge|z-E_n(\theta+p\cdot\om)|\ge\delta_n.
\end{equation}
The choice of $\mathrm C_{5}$ in Subsection~\ref{sec:induction-proof}
and \eqref{2026100908} give
\[
 \frac{\mathrm C_{5} a^2\beta^{-1}N_n^{\tau-2r}}{\delta_n}
 =\mathrm C_{5} a^2\beta^{-2}N_n^{2\tau-2r}\le\frac1{16}.
\]
Combining \eqref{2026101403} and \eqref{2026101404}
proves \eqref{2026101306}.
For $n=1$, \eqref{2026101403} proves the estimate
\eqref{2026100510} required in Subsection~\ref{sec:next-equations}.
The ranges $1\le h\le n-1$ and $1\le k<n$ are then empty.
The membership $B_1(p)\in\Upsilon_0$ is automatic, and
Lemma~\ref{lem:first-region-comparison} gives the comparison
for $B_1(p)\subset\Omega$.

We now solve for $\left.v\right|_{I_n}$, where $I_n=Q_n\cap\Omega$ and
$P_n=S_n\cap\Omega$ as in \eqref{2026100707}.
The inclusion $S_{n-1}\subset\Res_{n-1}$ and \eqref{2026100507}
give $M_n$-separation of $S_{n-1}$.
Assertion \textup{(4)} at $i=n-1$ gives
\begin{equation}\label{2026101501}
 |\cA_\Omega^{(n-1)}(x,y)|\le2a\langle x-y\rangle^{-r},
 \qquad x,y\in S_{n-1}\cap\Omega,\quad x\ne y.
\end{equation}
To verify the assumptions of Lemma~\ref{lem:sparse-schur} for
$\cA_\Omega^{(n-1)}$, take $\fX=S_{n-1}\cap\Omega$,
$\fP=P_n$, $\fQ=I_n$, $M=M_n$, $a_*=2a$ and
$\Delta=\delta_n/2$. The entry bounds are
\eqref{2026101306} and \eqref{2026101501}, and
\[
 a_*\Delta^{-1}M^{-r}=4a\delta_n^{-1}M_n^{-r}
 \le(16\mathrm C_{5})^{-1}\le(2\mathrm C_{3})^{-1}
\]
by \eqref{2026100907}.
Lemma~\ref{lem:sparse-schur} proves that
$(\cA_\Omega^{(n-1)})_{I_n}$ is invertible.
Thus $\cA_\Omega^{(n-1)}\left.v\right|_{P_{n-1}}=0$ is equivalent to
\begin{equation*}
 \left.v\right|_{I_n}=-(\cA_\Omega^{(n-1)})_{I_n}^{-1}
                (\cA_\Omega^{(n-1)})_{I_n,P_n}\left.v\right|_{P_n},\qquad \cA_\Omega^{(n)}\left.v\right|_{P_n}=0.
\end{equation*}
It remains to prove
$|\cA_\Omega^{(n)}(p,q)|\le2a\langle p-q\rangle^{-r}$ for
$p,q\in S_n\cap\Omega$, $p\ne q$.
Since $\Omega\in\Upsilon_n\subset\Upsilon_k$ for $0\le k\le n$,
assertion \textup{(4)} verifies the entry assumptions of
Lemma~\ref{lem:sparse-schur} at $0\le k<n$.
At $k=n$, \eqref{2026101306} and
\eqref{2026101501} verify the diagonal and off-diagonal
assumptions, respectively.
With $\fX=S_{k-1}\cap\Omega$, $\fP=S_k\cap\Omega$,
$\fQ=Q_k\cap\Omega$, $M=M_k$, $a_*=2a$ and
$\Delta=\delta_k/2$, the separation follows from \eqref{2026100507} for $k\ge1$
and from $M_0=1$ for $k=0$.
Condition \eqref{2026100907} verifies the
smallness assumption. Summing \eqref{2026102101} over $0\le k\le n$, we obtain
\begin{align}
 &|\cA_\Omega^{(n)}(p,q)-\cA_\Omega^{(-1)}(p,q)|
 \le\sum_{k=0}^{n}|\cA_\Omega^{(k)}(p,q)-\cA_\Omega^{(k-1)}(p,q)|\notag\\*
 &\le\mathrm C_{5} a^2\sum_{k=0}^{n}\delta_k^{-1}M_k^{-r}
                         \langle p-q\rangle^{-r}
 \le\frac a{64}\langle p-q\rangle^{-r}.
 \label{2026101502}
\end{align}
Combining \eqref{2026093004} and
\eqref{2026101502} proves
\eqref{2026101307} and extends assertion \textup{(4)} to $i=n$.
\end{proof}

\subsection{The next finite region and local eigenvalue}\label{sec:block-geometry}
We first construct $B_{n+1}$ satisfying
\eqref{2026101004}--\eqref{2026101007}, which provide the
inclusions and distances needed to compare the Schur complements.
For a finite $B\subset\mathbb Z^d$, a set $X\subset\mathbb Z^d$
and $R\ge1$, define
\begin{equation}\label{2026101503}
 \widetilde B=B\cup
    \bigcup_{\substack{p\in X\\(p+\Lambda_R)\cap B\ne\varnothing}}
            (p+\Lambda_R).
\end{equation}
Start with $B=\Lambda_{4N_{n+1}}$ and successively apply
\begin{equation}\label{2026101504}
 B\longmapsto\widetilde B,\qquad
 X=Z_h,\quad R=9N_{h+1},\quad h=n-1,n-2,\ldots,0.
\end{equation}
Denote the final set by $B_{n+1}$.
Since $B_{h+1}\subset\Lambda_{8N_{h+1}}$, the ball
$p+\Lambda_{9N_{h+1}}$ contains $B_{h+1}(p)+\Lambda_{N_{h+1}}$,
as required in \eqref{2026101005}.
The sets $Z_0,\ldots,Z_{n-1}$ and numbers $N_1,\ldots,N_{n+1}$
are fixed before $t,\theta,z$ are chosen.

\begin{proposition}\label{prop:det-blocks}
The set $B_{n+1}$ defined by \eqref{2026101504} is finite
and extends assertion \textup{(1)} of
Theorem~\ref{thm:finite-volume-induction} from step $n$ to step $n+1$.
\end{proposition}
The proof of Proposition~\ref{prop:det-blocks} is given in
Appendix~\ref{app:block-verification}.

We now work on $B_{n+1}$ and retain the component $v(o)$.
Define $J_{n+1}$ and $s_{n+1}$ by \eqref{2026100903} and
\eqref{2026100901} with $m=n+1$, provided that
$\cA_{B_{n+1}^\circ}(z;\theta)$ is invertible.

\begin{lemma}\label{lem:next-local-eigenvalue}
For every $\theta\in\Theta$ and $z\in J_{n+1}(\theta)$,
$\cA_{B_{n+1}^\circ}(z;\theta)$ is invertible.
The function $s_{n+1}(\cdot;\theta)$ has exactly one zero
$E_{n+1}(\theta)$ in $J_{n+1}(\theta)$ and
\begin{equation}\label{2026101505}
 |E_{n+1}(\theta)-E_n(\theta)|
 \le\mathrm C_{5} a^2\beta^{-1}N_n^{\tau-2r}\le\delta_n/16.
\end{equation}
The eigenvalue $E_{n+1}(\theta)$ is simple and satisfies
\eqref{2026101009} and
\eqref{2026101008} with $m=n+1$.
For the homotopy \eqref{2026100203}, the invertibility,
uniqueness and estimates in the lemma hold for every $t\in[0,1]$, and
$E_{n+1}(\theta;t)=\lambda(B_{n+1};\theta;t)$.
\end{lemma}
\begin{proof}
We first prove invertibility on $B_{n+1}^\circ$.
We then estimate $s_{n+1}(E_n(\theta);\theta)$ by comparing
$B_n$ with $B_{n+1}$ and eliminating $Q_n\cap B_{n+1}$.
Finally, Proposition~\ref{prop:homotopy-monotonicity} identifies
the new local eigenvalue for every $t\in[0,1]$.

Let $z\in J_{n+1}(\theta)$.
By \eqref{2026101102}, for $0\le h<n$,
\begin{equation}\label{2026101506}
 |z-E_h(\theta)|\le4\delta_n+\delta_h/8<\delta_h.
\end{equation}
The strict inequality follows from $\delta_n\le16^{-\tau}\delta_h$
and $\tau>d\ge1$.  Thus $o\in S_{n-1}$ throughout $J_{n+1}(\theta)$.
For $p\in\Res_h$ with $0\le h<n$, \eqref{2026100507} gives
\[
 L\|p\cdot\om\|_{\mathbb T}
 \le|E_h(\theta+p\cdot\om)-E_h(\theta)|
 <\delta_h+4\delta_n+\delta_h/8<8\delta_h.
\]
For $p\in\Res_n$, \eqref{2026100507} and
$|z-E_n(\theta)|\le4\delta_n$ give
\[
 L\|p\cdot\om\|_{\mathbb T}
 \le|E_n(\theta+p\cdot\om)-E_n(\theta)|<5\delta_n<8\delta_n.
\]
Hence $\Res_h\subset Z_h$ for $0\le h\le n$ and $z\in J_{n+1}(\theta)$.
The geometric inclusion \eqref{2026101005} and the exclusion
\eqref{2026101007} imply, respectively,
\begin{equation}\label{2026101507}
 B_{n+1}\in\Upsilon_n,\qquad B_{n+1}^\circ\in\Upsilon_n.
\end{equation}
Indeed, for $p\in B_{n+1}^\circ\cap S_{h-1}$ and $1\le h\le n$,
we have $p\in Z_{h-1}^\circ$, so
$B_h(p)\subset B_{n+1}$ and $o\notin B_h(p)$.
Also, \eqref{2026101006} gives
$B_{n+1}^\circ\cap S_n=\varnothing$.
Therefore
\[
 B_{n+1}^\circ=\bigsqcup_{h=0}^{n}(Q_h\cap B_{n+1}^\circ).
\]
By assertion \textup{(4)} and Lemma~\ref{lem:next-elimination},
$(\cA_{B_{n+1}^\circ}^{(h-1)})_{Q_h\cap B_{n+1}^\circ}$ is invertible
for $0\le h\le n$.
Since $B_{n+1}^\circ=\bigsqcup_{h=0}^n(Q_h\cap B_{n+1}^\circ)$,
back substitution for $h=n,n-1,\ldots,0$ in
$\cA_{B_{n+1}^\circ}(z;\theta)v=0$ gives $v=0$.
Thus $\cA_{B_{n+1}^\circ}(z;\theta)$ is invertible for every
$z\in J_{n+1}(\theta)$.
The equation on $B_{n+1}$ is equivalent to
\begin{equation*}
 \left.v\right|_{B_{n+1}^\circ}=-\cA_{B_{n+1}^\circ}(z;\theta)^{-1}c_{n+1}v(o),\qquad s_{n+1}(z;\theta)v(o)=0.
\end{equation*}
By Lemma~\ref{lem:schur-energy-comparison},
$s_{n+1}(z_2;\theta)-s_{n+1}(z_1;\theta)\ge z_2-z_1$
for $z_1<z_2$ in $J_{n+1}(\theta)$.
When $o\notin S_n$, the definition \eqref{2026100901}
still retains $v(o)$: the elimination at step $n$ uses
$Q_n\cap B_{n+1}^\circ$ and leaves the index set $\{o\}$.

Take $z=E_n(\theta)$ and let
$\cM=\cA_{B_{n+1}}^{(n-1)}(E_n(\theta);\theta)$.
By \eqref{2026101506}, $o\in S_{n-1}$, and
\eqref{2026101507} gives $B_{n+1}\in\Upsilon_n$.
Since $s_n(E_n(\theta);\theta)=0$ by assertion \textup{(2)},
\eqref{2026101403} with $p=o$ and $\Omega=B_{n+1}$ yields
\begin{equation}\label{2026101508}
 |\cM(o,o)-s_n(E_n(\theta);\theta)|
 \le\frac{\mathrm C_{5}}{2}a^2\beta^{-1}N_n^{\tau-2r}.
\end{equation}
Now $o\in S_n$, since $z=E_n(\theta)$.
Also $B_{n+1}\cap S_n=\{o\}$ by \eqref{2026101006}.
Let $I_n=Q_n\cap B_{n+1}$, so $o\notin I_n$.
Solving for $\left.v\right|_{I_n}$ gives
\begin{equation}\label{2026101509}
 s_{n+1}(E_n(\theta);\theta)
 =\bigl[\cM(o,o)-s_n(E_n(\theta);\theta)\bigr]
       -\cM_{o,I_n}\cM_{I_n}^{-1}\cM_{I_n,o}.
\end{equation}
The first term in \eqref{2026101509} is the change
after replacing $B_n$ by $B_{n+1}$, bounded in
\eqref{2026101508}.
The second term comes from eliminating $I_n$ in $B_{n+1}$.
Lemma~\ref{lem:next-elimination}, \eqref{2026101501}
and \eqref{2026100907} verify the assumptions
of Lemma~\ref{lem:sparse-schur} for $\cM$ with
$\fX=\{o\}\sqcup I_n$, $\fP=\{o\}$, $\fQ=I_n$,
$M=M_n=N_{n-1}$, $a_*=2a$ and $\Delta=\delta_n/2$;
the set $\{o\}\sqcup I_n\subset S_{n-1}$ is $M_n$-separated.
Equation \eqref{2026102102} therefore implies
\begin{equation}
 |\cM_{o,I_n}\cM_{I_n}^{-1}\cM_{I_n,o}|
 \le16\mathrm C_{3}a^2\delta_n^{-1}M_n^{-2r} =16^{2r+1}\mathrm C_{3}a^2\beta^{-1}N_n^{\tau-2r}
 \le\frac{\mathrm C_{5}}{2}a^2\beta^{-1}N_n^{\tau-2r}.
 \label{2026101601}
\end{equation}
Substituting \eqref{2026101508} and
\eqref{2026101601} into
\eqref{2026101509}, and using
\eqref{2026100904}, we obtain
\[
 |s_{n+1}(E_n(\theta);\theta)|
 \le\mathrm C_{5} a^2\beta^{-1}N_n^{\tau-2r}.
\]
By \eqref{2026100908}, the right-hand side is at most
$\delta_n/16$.
By Lemma~\ref{lem:schur-energy-comparison},
\begin{equation*}
 s_{n+1}(E_n(\theta)-4\delta_n;\theta)\le-63\delta_n/16<0,\qquad s_{n+1}(E_n(\theta)+4\delta_n;\theta)\ge63\delta_n/16>0.
\end{equation*}
The construction \eqref{2026101504} is independent of $t$,
and replacing $\eps$ by $t\eps$ preserves the bounds from
\eqref{2026101506} through
\eqref{2026101601} for $t\in[0,1]$.
Restoring $t$ in the notation,
\eqref{2026100902} implies
\[
 E_n(\theta;t)\pm4\delta_n\notin
 \Spec\cH_{B_{n+1}}(\theta;t),\qquad t\in[0,1].
\]
At $t=0$, $E_n(\theta;0)=f(\theta)$ by
\eqref{2026101203}. Since
$B_{n+1}\cap Z_n=\{o\}$,
\[
 |f(\theta+p\cdot\om)-f(\theta)|
 \ge L\|p\cdot\om\|_{\mathbb T}>8\delta_n,
 \qquad p\in B_{n+1}^\circ.
\]
By Proposition~\ref{prop:homotopy-monotonicity} with
$F=B_{n+1}$, $h(\theta;t)=E_n(\theta;t)$ and $\eta=4\delta_n$,
the function $E_{n+1}(\theta;t)=\lambda(B_{n+1};\theta;t)$ is
continuous in $t$ and is the unique eigenvalue in
$J_{n+1}(\theta;t)$, counted with multiplicity.
Equation \eqref{2026100902} identifies
$E_{n+1}(\theta;t)$ with the unique zero of $s_{n+1}(\cdot;\theta;t)$
in $J_{n+1}(\theta;t)$.

Fixing $t$ once more, Lemma~\ref{lem:schur-energy-comparison} gives
\begin{equation*}
 |E_{n+1}(\theta)-E_n(\theta)|
 \le |s_{n+1}(E_n(\theta);\theta)| \le\mathrm C_{5} a^2\beta^{-1}N_n^{\tau-2r}\le\delta_n/16,
\end{equation*}
which is \eqref{2026101505}.
Lemma~\ref{lem:schur-energy-comparison} also proves
\eqref{2026101008} at scale $n+1$.
Since
$\Spec\cH_{B_{n+1}}(\theta)\cap J_{n+1}(\theta)=\{E_{n+1}(\theta)\}$,
we have
\[
 \dist\!\left(E_{n+1}(\theta),
 \Spec\cH_{B_{n+1}}(\theta)\setminus\{E_{n+1}(\theta)\}\right)
 \ge4\delta_n-|E_{n+1}(\theta)-E_n(\theta)|\ge3\delta_n.
\]
This proves \eqref{2026101009} for every $t\in[0,1]$.
\end{proof}
\subsection{Eigenvector estimates and increments}
\label{sec:finite-eigenvectors}
At $z=E_j(\theta;t)$, we recover the coordinates on $B_j^\circ$
from the value $v_{j,\theta}(t)(o)=1$. We then compare the
zero-extended eigenvectors at scales $j$ and $j+1$.

\begin{lemma}\label{lem:finite-eigenvector-bound}
Let $j\ge1$ and suppose that assertions \textup{(1)--(4)} of
Theorem~\ref{thm:finite-volume-induction} hold through step $j$
at a fixed $t\in[0,1]$.
Then there is a unique eigenvector $v_{j,\theta}(t)$ of
$\cH_{B_j}(\theta;t)$ corresponding to $E_j(\theta;t)$
and satisfying $v_{j,\theta}(t)(o)=1$.
Extended by zero outside $B_j$, the vector satisfies
\eqref{2026101201}, uniformly in $\theta$.
If assertions \textup{(1)--(4)} hold for every $t\in[0,1]$,
the bound is uniform in $t$ as well.
\end{lemma}
\begin{proof}
Fix $\theta$ and $t$, suppress $t$, and put $z=E_j(\theta)$.
For $0\le h<j$, \eqref{2026101102} gives $o\in S_h$.
For $p\in S_h$, \eqref{2026100507} gives
\[
 L\|p\cdot\om\|_{\mathbb T}
 \le |E_h(\theta+p\cdot\om)-E_h(\theta)|
 <\delta_h+\delta_h/8<8\delta_h.
\]
Thus $S_h(E_j(\theta);\theta)\subset Z_h$. By
\eqref{2026101005}, \eqref{2026101006} and
\eqref{2026100803}--\eqref{2026100804},
\begin{equation}\label{2026101602}
 B_j\in\Upsilon_{j-1},\qquad B_j\cap S_{j-1}=\{o\}.
\end{equation}
Use $I_i,P_i$ from \eqref{2026100707} with
$\Omega=B_j$. Then $o\in P_i$ for $-1\le i<j$ and
$P_{j-1}=\{o\}$.
The final equation is
$s_j(E_j(\theta);\theta)v(o)=0$, so set $v^{(j-1)}=e_o$.
Assertion \textup{(4)} permits the following back substitution,
for $i=j-1,\ldots,0$:
\begin{equation}\label{2026101603}
 v^{(i-1)}=\cW_{B_j}^{(i)}v^{(i)},\qquad
 \cW_{B_j}^{(i)}=
 \begin{pmatrix}
  \Id_{P_i}\\
  -(\cA_{B_j}^{(i-1)})_{I_i}^{-1}
       (\cA_{B_j}^{(i-1)})_{I_i,P_i}
 \end{pmatrix}
 :\mathbb C^{P_i}\longrightarrow\mathbb C^{P_{i-1}}.
\end{equation}
The row blocks are ordered as $P_{i-1}=P_i\sqcup I_i$.
Equations \eqref{2026100801} and
\eqref{2026100802} give
\[
 \cA_{B_j}^{(i)}v^{(i)}=0
 \quad\Longrightarrow\quad
 \cA_{B_j}^{(i-1)}v^{(i-1)}=0.
\]
Since $o\in P_i$, the back substitution preserves $v^{(i)}(o)=1$.
Consequently, the zero extension of
$\cW_{B_j}^{(0)}\cdots\cW_{B_j}^{(j-1)}e_o$
satisfies \eqref{2026101106}.
The invertibility of $\cA_{B_j^\circ}(E_j(\theta);\theta)$
gives uniqueness with this normalization.

For $i\ge1$, \eqref{2026100507} shows that $P_{i-1}$ is
$M_i$-separated; for $i=0$, use $M_0=1$.
The bounds \eqref{2026101104}--\eqref{2026101105}
and \eqref{2026100907} verify the assumptions
of Lemma~\ref{lem:sparse-schur} for $\cA_{B_j}^{(i-1)}$ with
$\fX=P_{i-1}$, $\fP=P_i$, $\fQ=I_i$, $M=M_i$,
$a_*=2a$ and $\Delta=\delta_i/2$.
Since $\rho<r-\tau<r-d$, \eqref{2026102104} gives
\begin{align}
 &\matnorm{\cW_{B_j}^{(0)}}{\rho}
 \le1+\mathrm C_{4}a\beta^{-1},
 \label{2026101604}\\
 &\matnorm{\cW_{B_j}^{(i)}}{\rho}
 \le1+\mathrm C_{4}a\delta_i^{-1}N_{i-1}^{\rho-r}
 =1+\mathrm C_{4}a\beta^{-1}16^\tau N_{i-1}^{\tau+\rho-r},
 \quad 1\le i<j.
 \label{2026101605}
\end{align}
The equality $v_{j,\theta}(o)=1$ implies
$\|v_{j,\theta}-e_o\|_{1,\rho}=\|v_{j,\theta}\|_{1,\rho}-1$.
Combining \eqref{2026101604},
\eqref{2026101605} and \eqref{2026092902}, we obtain
\begin{equation}\label{2026101606}
 \|v_{j,\theta}-e_o\|_{1,\rho}
 \le\prod_{i=0}^{j-1}\matnorm{\cW_{B_j}^{(i)}}{\rho}-1
 \le\exp\bigl(\mathrm C_{4}a\beta^{-1}G(\rho)\bigr)-1
 \le\mathrm C_{6} a.
\end{equation}
Here \eqref{2026100906} and
\eqref{2026100909} give
$\mathrm C_{4}a\beta^{-1}G(\rho)\le1/8$, so
$\exp(x)-1\le2x$ for $0\le x\le1/8$ proves the last inequality.
All bounds are independent of $\theta$ and $t$.
\end{proof}

\begin{proposition}\label{prop:finite-vector-increment}
Let $j\ge1$ and suppose that assertions \textup{(1)--(4)} of
Theorem~\ref{thm:finite-volume-induction} hold through step $j+1$
for every $t\in[0,1]$, with constants independent of $t$.
The eigenvectors constructed in Lemma~\ref{lem:finite-eigenvector-bound}
satisfy
\[
 \sup_{\theta\in\Theta,\,t\in[0,1]}
 \|v_{j+1,\theta}(t)-v_{j,\theta}(t)\|_1
 \le12\mathrm C_{6}aN_j^{-\kappa},
\]
where $\mathrm C_{6}$ and $\kappa$ are fixed in
\eqref{2026100906} and
\eqref{2026101002}.
\end{proposition}
\begin{proof}
Fix $\theta\in\Theta$ and $t\in[0,1]$, and suppress $t$ in
$v_{i,\theta}$ and $E_i$.
We compare $v_{j+1,\theta}$ with $v_{j,\theta}$ by expanding
$v_{j+1,\theta}$ in the translates
$v_{j,\theta+k\cdot\om}(\cdot-k)$, $k\in\mathbb Z^d$.
For $i=j,j+1$ and $m,k\in\mathbb Z^d$, define
\begin{align}
 &\mathcal V^{(i)}(m,k)=v_{i,\theta+k\cdot\om}(m-k),\qquad \mathcal D^{(i)}=\operatorname{diag}_{k\in\mathbb Z^d}
                         E_i(\theta+k\cdot\om),
 \label{2026101607}\\
 &\mathcal R^{(i)}(m,k)
 =t\eps\,\mathbf1_{\{m-k\notin B_i\}}
       \sum_{\ell\in\mathbb Z^d}\phi(m-\ell)\mathcal V^{(i)}(\ell,k).
 \label{2026101608}
\end{align}
By \eqref{2026092708} and
\eqref{2026101106},
\begin{equation}\label{2026101701}
 \cH(\theta;t)\mathcal V^{(i)}
 =\mathcal V^{(i)}\mathcal D^{(i)}+\mathcal R^{(i)}
\end{equation}
on finitely supported vectors.
Lemma~\ref{lem:finite-eigenvector-bound},
\eqref{2026092901} and Lemma~\ref{lem:column-algebra} give
\begin{equation}
 \matnorm{\mathcal V^{(i)}-\Id}{\rho}\le\mathrm C_{6} a\le1/4,\qquad \matnorm{(\mathcal V^{(i)})^{-1}}{\rho}
 \le(1-\mathrm C_{6} a)^{-1}\le4/3.
 \label{2026101702}
\end{equation}

Set $x=(\mathcal V^{(j)})^{-1}v_{j+1,\theta}$.
Since $\mathcal V^{(j)}e_o=v_{j,\theta}$, we have
\begin{equation}\label{2026101703}
 v_{j+1,\theta}-v_{j,\theta}=\mathcal V^{(j)}(x-e_o).
\end{equation}
For a radius $1\le R\le N_j$ to be chosen below, eigenvalue
separation will bound $\sum_{0<|k|\le R}|x(k)|$, while the
weighted norm will bound $\sum_{|k|>R}|x(k)|$.
We first estimate $\mathcal R^{(i)}$ and verify that
\eqref{2026101701} holds for $x$.
By translation covariance, it suffices to consider the column $k=o$.
The inclusion $\Lambda_{4N_i}\subset B_i$ implies
\[
 \sum_{m\notin B_i}\langle m-\ell\rangle^{-r}
 \le\mathrm C_{3}N_i^{d-r},\qquad |\ell|\le2N_i.
\]
For the remaining summation indices,
\eqref{2026101201} gives
\[
 \sum_{|\ell|>2N_i}|v_{i,\theta}(\ell)|
 \le\mathrm C_{6} a(1+2N_i)^{-\rho}.
\]
Splitting the sum in \eqref{2026101608}
at $|\ell|=2N_i$ and using
$\sum_m\langle m\rangle^{-r}\le\mathrm C_{3}$ and $\rho<r-d$, we obtain
\begin{equation}
 \matnorm{\mathcal R^{(i)}}{0}
 \le\mathrm C_{3}a\bigl((1+\mathrm C_{6} a)N_i^{d-r}
                         +\mathrm C_{6} aN_i^{-\rho}\bigr) \le2\mathrm C_{3}aN_i^{-\rho}.
 \label{2026101704}
\end{equation}

To use \eqref{2026101701} for $x$, set
\begin{equation}\label{2026101705}
 \mathcal F=\operatorname{diag}_{k\in\mathbb Z^d}f(\theta+k\cdot\om)
\end{equation}
on $\ell^1$, with the natural domain.
By \eqref{2026101102}, $\mathcal D^{(i)}-\mathcal F$ is bounded
on $\ell^1$, so
$\operatorname{Dom}(\mathcal D^{(i)})=\operatorname{Dom}(\mathcal F)$.
The identity \eqref{2026101701} extends
to $\operatorname{Dom}(\mathcal F)$, and
\begin{equation}\label{2026101706}
 \mathcal V^{(i)}\operatorname{Dom}(\mathcal F)
 =\operatorname{Dom}(\mathcal F),\qquad i=j,j+1.
\end{equation}
The domain assertion and the extension of
\eqref{2026101701} are proved in
Appendix~\ref{app:domain-preservation}.
By \eqref{2026101702},
\eqref{2026101201} and \eqref{2026100909},
\begin{equation}
 \|x\|_1\le\frac{1+\mathrm C_{6} a}{1-\mathrm C_{6} a},\qquad \|x-e_o\|_{1,\rho}
 \le\frac{2\mathrm C_{6} a}{1-\mathrm C_{6} a}
 \le3\mathrm C_{6}a.
 \label{2026101707}
\end{equation}
The vector $v_{j+1,\theta}$ has finite support; hence
\eqref{2026101706} implies
$x\in\operatorname{Dom}(\mathcal F)$.
Using \eqref{2026101701} for $i=j,j+1$, we obtain
\begin{equation}\label{2026101801}
 (\mathcal D^{(j)}-E_{j+1}(\theta))x
 =(\mathcal V^{(j)})^{-1}
       \bigl(\mathcal R^{(j+1)}e_o-\mathcal R^{(j)}x\bigr).
\end{equation}
Equations \eqref{2026101702},
\eqref{2026101704} and
\eqref{2026101707}, together with $N_{j+1}>N_j$,
give
\[
 \bigl\|(\mathcal V^{(j)})^{-1}
       (\mathcal R^{(j+1)}e_o-\mathcal R^{(j)}x)\bigr\|_1
 \le8\mathrm C_{3}aN_j^{-\rho}.
\]

For $1\le R\le N_j$ and $0<|k|\le R$,
\eqref{2026100507}, \eqref{2026101101},
\eqref{2026100908} and \eqref{2026093001} give
\begin{equation*}
 |E_j(\theta+k\cdot\om)-E_{j+1}(\theta)|
 \ge L\gamma\langle k\rangle^{-\tau}-\delta_j/16 \ge\frac{L\gamma}{2}(1+R)^{-\tau}.
\end{equation*}
Consequently, \eqref{2026101801} implies
\[
 \sum_{0<|k|\le R}|x(k)|
 \le\frac{2^{\tau+4}\mathrm C_{3}}{L\gamma}
             aR^\tau N_j^{-\rho}
 \le\mathrm C_{6} aR^\tau N_j^{-\rho}.
\]
For $|k|>R$, \eqref{2026101707} gives
\[
 \sum_{|k|>R}|x(k)|\le3\mathrm C_{6}aR^{-\rho}.
\]
Adding the estimates for $0<|k|\le R$ and $|k|>R$, we obtain
\begin{equation}\label{2026101802}
 \sum_{k\ne o}|x(k)|
 \le3\mathrm C_{6} a\bigl(R^\tau N_j^{-\rho}+R^{-\rho}\bigr).
\end{equation}
Taking $R=N_j^{\rho/(\tau+\rho)}\in[1,N_j]$ in
\eqref{2026101802}, we obtain
$\sum_{k\ne o}|x(k)|\le6\mathrm C_{6} aN_j^{-\kappa}$.

Finally, $v_{j+1,\theta}(o)=v_{j,\theta}(o)=1$ and
$v_{j+1,\theta}=\mathcal V^{(j)}x$ imply
\[
 x(o)-1=-\sum_{k\ne o}\mathcal V^{(j)}(o,k)x(k),\qquad
 |x(o)-1|\le\mathrm C_{6} a\sum_{k\ne o}|x(k)|.
\]
Using \eqref{2026101703}, we conclude that
\begin{equation*}
 \|v_{j+1,\theta}-v_{j,\theta}\|_1
 \le(1+\mathrm C_{6} a)^2\sum_{k\ne o}|x(k)| \le12\mathrm C_{6}aN_j^{-\kappa}.
\end{equation*}
The estimates are uniform in $\theta$ and $t$ and apply to $j=1$.
\end{proof}

\subsection{Completion of the induction}

\begin{proposition}\label{prop:general-inductive-step}
Let $n\ge1$. Suppose that $B_m$, $E_m$ and $v_{m,\theta}$,
$0\le m\le n$, satisfy assertions
\textup{(1)--(6)} of Theorem~\ref{thm:finite-volume-induction}
at step $n$, with the fixed constants.
Then there exist $B_{n+1}$, $E_{n+1}$ and $v_{n+1,\theta}$
such that assertions \textup{(1)--(6)} of
Theorem~\ref{thm:finite-volume-induction} hold at step $n+1$,
with the constants fixed in Subsection~\ref{sec:induction-proof}
and without changing $(B_m,E_m,v_{m,\theta})$ for $0\le m\le n$.
\end{proposition}
\begin{proof}
Proposition~\ref{prop:det-blocks} proves assertion \textup{(1)} for $B_{n+1}$.
The proof of Lemma~\ref{lem:next-elimination} applies for every
$t\in[0,1]$ and proves assertion \textup{(4)} at $i=n$.
By Lemma~\ref{lem:next-local-eigenvalue}, assertion \textup{(2)}
and \eqref{2026101101} hold at the new step for every $t\in[0,1]$.

By \eqref{2026101003} and Proposition~\ref{prop:kps-label},
\[
 E_{n+1}(\theta_2;t)-E_{n+1}(\theta_1;t)\ge L(\theta_2-\theta_1),
 \qquad \theta_1<\theta_2,\quad \theta_1,\theta_2\in\Theta\cap[0,1),\quad t\in[0,1].
\]
This proves \eqref{2026101103} at the new step.

For $0\le\ell<m\le n+1$, summing \eqref{2026101101} and using
\eqref{2026100908} and $N_\ell^{2\tau-2r}\le1$, we obtain
\begin{equation*}
 \|E_m(\cdot;t)-E_\ell(\cdot;t)\|_\infty
 \le\mathrm C_{5} a^2\beta^{-1}\sum_{h=\ell}^{m-1}N_h^{\tau-2r} \le\frac{\mathrm C_{5} a^2\beta^{-1}N_\ell^{\tau-2r}}
                  {1-16^{\tau-2r}}\le\delta_\ell/16.
\end{equation*}
This proves \eqref{2026101102} and completes assertion \textup{(3)}.
By Lemma~\ref{lem:finite-eigenvector-bound} with $j=n+1$,
there is an eigenvector $v_{n+1,\theta}(t)$ satisfying
\eqref{2026101201}.
By Proposition~\ref{prop:finite-vector-increment} with $j=n$,
$v_{n+1,\theta}(t)$ also satisfies \eqref{2026101202},
which proves assertion \textup{(5)}.
The normalization is $v_{n+1,\theta}(t)(o)=1$, and
\[
 \left.v_{n+1,\theta}(t)\right|_{B_{n+1}^\circ}
 =-\cA_{B_{n+1}^\circ}(E_{n+1}(\theta;t);\theta;t)^{-1}c_{n+1}(t).
\]
Assertion \textup{(2)} gives the invertibility of
$\cA_{B_{n+1}^\circ}(E_{n+1}(\theta;t);\theta;t)$ for every $t\in[0,1]$.
The identity $E_{n+1}(\theta;t)=\lambda(B_{n+1};\theta;t)$ from
Lemma~\ref{lem:next-local-eigenvalue}, together with
\eqref{2026100404} and \eqref{2026092905}, gives
$E_{n+1}(\theta;0)=f(\theta)$ and continuity in $t$.
The matrix $\cA_{B_{n+1}^\circ}(E_{n+1}(\theta;t);\theta;t)$ and
$c_{n+1}(t)$ are therefore continuous in $t$, and continuity of matrix inversion proves
continuity of $v_{n+1,\theta}(t)$.
Moreover, $c_{n+1}(0)=0$ gives $v_{n+1,\theta}(0)=e_o$.
The set $B_{n+1}$ and the constants in
Subsection~\ref{sec:induction-proof} are independent of $t$,
so assertion \textup{(6)} holds at step $n+1$.

\end{proof}

\begin{proof}[Proof of Theorem~\ref{thm:finite-volume-induction}]
The verification at the end of Section~\ref{sec:inductive-hypotheses}
establishes all six assertions for $n=1$.
Proposition~\ref{prop:general-inductive-step} extends
$(B_m,E_m,v_{m,\theta})_{0\le m\le n}$ from each $n\ge1$ to $n+1$.
Induction therefore constructs the eigenvalues $E_m(\cdot;t)$
and corresponding eigenvectors $v_{m,\theta}(t)$
for all $m\ge0$ and $t\in[0,1]$.
The common choice \eqref{2026101001} guarantees all the
smallness assumptions, including \eqref{2026100909}.
\end{proof}

\section{Localization}\label{sec:localization}
We use the summable increments in
Theorem~\ref{thm:finite-volume-induction} to construct eigenvectors
of $\cH(\theta)$.
We then prove that the normalized translates of $v_{\theta+n\cdot\om}$,
$n\in\mathbb Z^d$, form an orthonormal basis and establish
\eqref{2026092607}.
Throughout this section, use $\rho$ and $\eps_0$ fixed in
Subsection~\ref{sec:induction-proof}, assume $|\eps|<\eps_0$ and take $t=1$.
We omit the parameter $t=1$ from the notation.

\subsection{Limiting eigenvectors}\label{sec:limiting-eigenvectors}
By \eqref{2026101101},
\[
 \sum_{j\ge1}\|E_j-E_{j-1}\|_\infty
 \le\mathrm C_{5}a^2\beta^{-1}
          \sum_{j\ge1}N_{j-1}^{\tau-2r}<\infty.
\]
Since $E_0=f$, the functions $E_j-f$ converge uniformly on $\fT$.
We define
\[
 E_\infty=f+\lim_{j\to\infty}(E_j-f).
\]
The function $E_\infty:\fT\to\mathbb R$ is $1$-periodic.
Summing \eqref{2026101101} over $j\ge1$, we obtain
\[
 \|E_\infty-f\|_\infty
 \le\frac{\mathrm C_{5}a^2}{\beta(1-16^{\tau-2r})}
 \le\mathrm C_{1}|\eps|^2.
\]
Here $\mathrm C_{1}$ depends only on $d,r,\tau,\gamma,L,\phi$
and is independent of $\rho$, as required in Theorem~\ref{thm:main}.
By \eqref{2026101103},
\[
 E_\infty(\theta_2)-E_\infty(\theta_1)
 =\lim_{j\to\infty}(E_j(\theta_2)-E_j(\theta_1))\ge L(\theta_2-\theta_1),
 \quad \theta_1,\theta_2\in\fT\cap[0,1),\quad \theta_1<\theta_2.
\]
By \eqref{2026100507} and \eqref{2026092502},
we obtain
\begin{equation}\label{2026101803}
 |E_\infty(\theta+n\cdot\om)-E_\infty(\theta+m\cdot\om)|
 \ge L\|(n-m)\cdot\om\|_{\mathbb T}
 \ge L\gamma\langle n-m\rangle^{-\tau},
 \qquad n\ne m.
\end{equation}
We now use \eqref{2026101201} and
\eqref{2026101202} to construct an
eigenvector of $\cH(\theta)$ with eigenvalue $E_\infty(\theta)$.

\begin{lemma}\label{lem:limiting-eigenvector}
For every $\theta\in\fT$, there is a vector
$v_\theta\in\operatorname{Dom}(\cH(\theta))$ such that
\begin{equation*}
 \cH(\theta) v_\theta=E_\infty(\theta)v_\theta,\qquad
 v_\theta({o})=1,\qquad \sum_{m\ne{o}}|v_\theta(m)|\langle m\rangle^\rho
 \le\mathrm C_{6}a\le\frac14.
\end{equation*}
\end{lemma}
\begin{proof}
Fix $\theta\in\fT$ and use the vectors in assertion
\textup{(\ref{item:induction-vectors})} of
Theorem~\ref{thm:finite-volume-induction}.
Since $N_j=16^j$ and $\kappa>0$, \eqref{2026101202} gives
\[
 \sum_{j\ge1}\sup_{\theta\in\Theta}
       \|v_{j+1,\theta}-v_{j,\theta}\|_1<\infty.
\]
Thus $v_{j,\theta}$ converges in $\ell^1$, uniformly in $\theta$,
to a vector $v_\theta$ with $v_\theta(o)=1$, and
\[
 \sup_{\theta\in\Theta}\|v_\theta-v_{j,\theta}\|_1
 \le\frac{12\mathrm C_{6}}{1-16^{-\kappa}}aN_j^{-\kappa}.
\]
In particular, the convergence also holds in $\ell^2$.
By Fatou's lemma, \eqref{2026101201} and
\eqref{2026100909},
\[
 \sum_{m\ne o}|v_\theta(m)|\langle m\rangle^\rho
 \le\liminf_{j\to\infty}
       \sum_{m\ne o}|v_{j,\theta}(m)|\langle m\rangle^\rho
 \le\mathrm C_{6}a\le\frac14.
\]

To verify the eigenvalue equation, fix $m\in\mathbb Z^d$.
By \eqref{2026101004}, $m\in\fB_j$ for all sufficiently
large $j$.  Equation \eqref{2026101106} gives
\[
 f(\theta+m\cdot\om)v_{j,\theta}(m)
 +\eps(\cT v_{j,\theta})(m)
 =E_j(\theta)v_{j,\theta}(m).
\]
Since $\theta\in\Theta$, the number $f(\theta+m\cdot\om)$ is finite.
The operator $\cT$ defined in \eqref{2026092707}
is bounded on $\ell^2$ by \eqref{2026092601}.
Hence $v_{j,\theta}\to v_\theta$ in $\ell^2$ and
$E_j(\theta)\to E_\infty(\theta)$ imply
\[
 f(\theta+m\cdot\om)v_\theta(m)
 =E_\infty(\theta)v_\theta(m)
  -\eps(\cT v_\theta)(m).
\]
Since
$E_\infty(\theta)v_\theta-\eps\cT v_\theta\in\ell^2(\mathbb Z^d)$,
we have $v_\theta\in\operatorname{Dom}(\cH(\theta))$.
Thus $\cH(\theta) v_\theta=E_\infty(\theta)v_\theta$.
\end{proof}

\subsection{Completeness and dynamical localization}\label{sec:completeness-dynamics}
We prove completeness by estimating the matrix of normalized
eigenvectors in the weighted column norm.

\begin{proof}[Proof of Theorem~\ref{thm:main}]
Use $\eps_0$ from \eqref{2026101001} and fix $\theta\in\fT$.
For each $n\in\mathbb Z^d$, let $v_{\theta+n\cdot\om}$ be the
eigenvector supplied by Lemma~\ref{lem:limiting-eigenvector} at
phase $\theta+n\cdot\om$, and define
\begingroup\postdisplaypenalty=100
\[
 \psi_{\theta,n}
 =\frac{\cU(-n)v_{\theta+n\cdot\om}}
        {\|v_{\theta+n\cdot\om}\|_2}.
\]
\endgroup
Here $\cU(-n)$ is the translation defined in
\eqref{2026092706}.
By \eqref{2026092708}, $\psi_{\theta,n}$ is a normalized
eigenvector of $\cH(\theta)$ with eigenvalue
$E_\infty(\theta+n\cdot\om)$ and $\psi_{\theta,n}(n)>0$.
By \eqref{2026101803},
$E_\infty(\theta+n\cdot\om)\ne E_\infty(\theta+m\cdot\om)$ for
$n\ne m$.  Self-adjointness of $\cH(\theta)$ therefore implies
\[
 \langle\psi_{\theta,n},\psi_{\theta,m}\rangle=\delta_{nm},
 \qquad n,m\in\mathbb Z^d.
\]

We next estimate the change caused by normalization.
Since $v_{\theta+n\cdot\om}({o})=1$, we have
\begin{equation}
 0\le1-\|v_{\theta+n\cdot\om}\|_2^{-1}
 \le\tfrac12\bigl(\|v_{\theta+n\cdot\om}\|_2^2-1\bigr) \le\tfrac12\left(\sum_{m\ne{o}}
        |v_{\theta+n\cdot\om}(m)|\right)^2
 \le\frac{\mathrm C_{6}^2}{2}a^2.
 \label{2026101804}
\end{equation}
Since $\psi_{\theta,n}(n)=\|v_{\theta+n\cdot\om}\|_2^{-1}$ and
$\sum_{m\ne o}|v_{\theta+n\cdot\om}(m)|\le1/4$,
\eqref{2026101804} and
Lemma~\ref{lem:limiting-eigenvector} imply
\begin{equation}\label{2026101805}
 \sup_n\sum_m|\psi_{\theta,n}(m)-\delta_{mn}|
                  \langle m-n\rangle^\rho
 \le2\mathrm C_{6}a.
\end{equation}
By \eqref{2026101805}, \eqref{2026093003}
and \eqref{2026092602}, we obtain \eqref{2026092606}
with a finite constant $\mathrm C_{2}$ depending only on
$d,r,\tau,\gamma,L,\rho$ and $\phi$.

Define the matrix $\Psi(\theta)$ by
\begin{equation}\label{2026101806}
 \Psi(\theta)(m,n)=\psi_{\theta,n}(m),\qquad
 \Psi(\theta)e_n=\psi_{\theta,n}.
\end{equation}
For a finitely supported vector $x$, orthonormality gives
$\|\Psi(\theta)x\|_2=\|x\|_2$.
Hence \eqref{2026101806} extends to an isometry on
$\ell^2(\mathbb Z^d)$.
By Lemma~\ref{lem:limiting-eigenvector},
\eqref{2026101804} and \eqref{2026092901},
we obtain
\begin{equation}\label{2026101901}
 \matnorm{\Psi(\theta)-\Id}{\rho}
 \le\frac14+\frac12\left(\frac14\right)^2
 =\frac9{32}<\frac12.
\end{equation}
By Lemma~\ref{lem:column-algebra}, $\Psi(\theta)$ is then invertible
on $\ell^1$, with
\[
 \Psi(\theta)^{-1}=\sum_{k=0}^\infty(\Id-\Psi(\theta))^k,
\]
where the series converges in the weighted column norm.
Moreover,
\begin{equation}\label{2026101902}
 \matnorm{\Psi(\theta)^{-1}}{\rho}
 \le\bigl(1-\matnorm{\Psi(\theta)-\Id}{\rho}\bigr)^{-1}\le2.
\end{equation}

The actions of $\Psi(\theta)$ on $\ell^1$ and $\ell^2$ agree
on $\ell^1$: approximate an $\ell^1$ vector by finitely supported
vectors and use the continuous inclusion $\ell^1\subset\ell^2$.

For every $n\in\mathbb Z^d$, we therefore have
$\Psi(\theta)^{-1}e_n\in\ell^1\subset\ell^2$ and
$\Psi(\theta)(\Psi(\theta)^{-1}e_n)=e_n$ in $\ell^2$.
The range of the isometry $\Psi(\theta)$ is closed and contains
every $e_n$, so $\operatorname{Ran}\Psi(\theta)=\ell^2(\mathbb Z^d)$.
Hence $\Psi(\theta)$ is unitary and
$(\psi_{\theta,n})_{n\in\mathbb Z^d}$ is a complete orthonormal basis.  In particular, $\cH(\theta)$ has pure
point spectrum.
For every $x\in\ell^1$, the identity
$\Psi(\theta)\Psi(\theta)^{-1} x=x$ in $\ell^2$ now implies
$\Psi(\theta)^{-1} x=\Psi(\theta)^*x$.

To prove \eqref{2026092607}, the spectral theorem gives, for $u\in\mathbb R$,
\begin{equation}\label{2026101903}
 \exp(-\mathrm i u\cH(\theta))
 =\Psi(\theta)\operatorname{diag}_{n\in\mathbb Z^d}
   \bigl(\exp(-\mathrm i uE_\infty(\theta+n\cdot\om))\bigr)\Psi(\theta)^*.
\end{equation}
For every $u\in\mathbb R$,
\[
 \matnorm{\operatorname{diag}_{n\in\mathbb Z^d}
   (\exp(-\mathrm i uE_\infty(\theta+n\cdot\om)))}{\rho}=1.
\]
On $\ell^1$, \eqref{2026101903}
therefore holds with $\Psi(\theta)^{-1}$ in place of $\Psi(\theta)^*$.
By
\eqref{2026092902}, \eqref{2026101901}
and \eqref{2026101902}, we obtain
\begin{equation*}
 \matnorm{\exp\left(-\mathrm i u\cH(\theta)\right)}{\rho}
 \le\matnorm{\Psi(\theta)}{\rho}\matnorm{\Psi(\theta)^{-1}}{\rho} \le\frac32\cdot2=3.
\end{equation*}
This proves \eqref{2026092607} and completes the proof.
\end{proof}

\appendix
\section{Sparse sums and Schur estimates}\label{app:matrix-estimates}
We estimate sums over $M$-separated sets and use the resulting
bounds to control one Schur elimination.

\subsection{Sparse sums}

\begin{lemma}[Sparse sums]\label{lem:sparse-sums}
Let $M\ge1$, let $\fX\subset\mathbb Z^d$ be $M$-separated and let
$r>d$. Use $\mathrm C_{3}(d,r)$ from \eqref{2026093005}.
For every $s>d$ and $p\in\fX$, we have
\begin{equation}\label{2026101904}
 \sum_{x\in\fX\setminus\{p\}}\langle p-x\rangle^{-s}
 \le\frac{5^d}{1-2^{d-s}}M^{-s}.
\end{equation}
For distinct $p,q\in\fX$, we have
\begin{equation}\label{2026101905}
 \sum_{x\in\fX\setminus\{p,q\}}
       \langle p-x\rangle^{-r}\langle x-q\rangle^{-r}
 \le\mathrm C_{3}(d,r)M^{-r}\langle p-q\rangle^{-r}.
\end{equation}
If $p\in\fX$, $\fY\subset\fX\setminus\{p\}$, $y\in\mathbb R^d$ and
$\dist(y,\fX)\ge M$, then
\begin{equation}\label{2026101906}
 \sum_{x\in\fY}\langle p-x\rangle^{-r}\langle x-y\rangle^{-r}
 \le\mathrm C_{3}(d,r)M^{-r}\langle p-y\rangle^{-r}.
\end{equation}
For every $p\in\mathbb Z^d$ and $R\ge M$, we also have
\begin{equation}\label{2026101907}
 \sum_{\substack{x\in\fX\\|x-p|\ge R}}\langle p-x\rangle^{-2r}
 \le\mathrm C_{3}M^{-d}R^{d-2r}.
\end{equation}
\end{lemma}
\begin{proof}
The half-open cubes $x+[-M/2,M/2)^d$, $x\in\fX$, are pairwise disjoint.
Comparison of volumes gives, for $p\in\mathbb R^d$ and $R\ge M$,
\begin{equation}\label{2026102001}
 \#\{x\in\fX:|x-p|<R\}
 \le\left(\frac{2R}{M}+1\right)^d
 \le3^d(R/M)^d.
\end{equation}
For $j\ge0$, \eqref{2026102001} implies
\[
 \#\{x\in\fX:2^jR\le|x-p|<2^{j+1}R\}
 \le\left(2^{j+2}\frac RM+1\right)^d
 \le5^d2^{jd}(R/M)^d.
\]
Consequently, for every $s>d$,
\begin{equation}
 \sum_{\substack{x\in\fX\\|x-p|\ge R}}\langle p-x\rangle^{-s}
 \le5^d(R/M)^d\sum_{j\ge0}2^{jd}(2^jR)^{-s} =\frac{5^d}{1-2^{d-s}}M^{-d}R^{d-s}.
 \label{2026102002}
\end{equation}
Taking $s=2r$ gives the tail bound \eqref{2026101907}.
For $p\in\fX$, all points of $\fX\setminus\{p\}$ have distance at least
$M$ from $p$, so \eqref{2026102002} with $R=M$ proves
\eqref{2026101904}.
Equation \eqref{2026102002} with $p=y$, $R=M$ and $s=r$ gives
\begin{equation}\label{2026102003}
 \sum_{x\in\fX}\langle x-y\rangle^{-r}
 \le\frac{5^d}{1-2^{d-r}}M^{-r},
 \qquad \dist(y,\fX)\ge M.
\end{equation}

For any $p,x,q\in\mathbb R^d$, the triangle inequality implies
$|p-x|\ge|p-q|/2$ or $|x-q|\ge|p-q|/2$.  Hence
\[
 \langle p-x\rangle^{-r}\langle x-q\rangle^{-r}
 \le2^r\langle p-q\rangle^{-r}
 \bigl(\langle p-x\rangle^{-r}+\langle x-q\rangle^{-r}\bigr).
\]
For distinct $p,q\in\fX$, summing over
$x\in\fX\setminus\{p,q\}$ and applying \eqref{2026101904} with
$s=r$ proves \eqref{2026101905}, with coefficient
$2^{r+1}5^d/(1-2^{d-r})$.
Replacing $q$ by $y$ and summing over $x\in\fY$, we obtain
\eqref{2026101906} from \eqref{2026101904} and
\eqref{2026102003}. Also, \eqref{2026101904} with
$\fX=\mathbb Z^d$, $M=1$, $s=r$ and $p=o$ bounds
$\sum_k\langle k\rangle^{-r}$ by $1+5^d/(1-2^{d-r})$.
The convolution, tail and full-lattice coefficients are finite
and depend only on $d,r$, which justifies the common choice of
$\mathrm C_3$ in \eqref{2026093005}.
\end{proof}

\subsection{One Schur elimination}
Use $\mathrm C_{3}$ from \eqref{2026093005} and the coefficients
$\mathrm C_{4}(d,r,s)$ fixed before \eqref{2026100307}.

\begin{lemma}\label{lem:sparse-schur}
Let $r>d$ and let $\fX=\fP\sqcup \fQ\subset\mathbb Z^d$ be finite and $M$-separated,
where $M\ge1$.  Let $\cM$ be a Hermitian matrix on $\fX$ such that
\[
 |\cM(p,q)|\le a_*\langle p-q\rangle^{-r}\quad(p\ne q),\qquad
 |\cM(q,q)|\ge\Delta\quad(q\in \fQ),
\]
where $a_*\ge0$ and $\Delta>0$.
If
\begin{equation}\label{2026102004}
 a_*\Delta^{-1}M^{-r}\le(2\mathrm C_{3})^{-1},
\end{equation}
then $\cM_\fQ$ is invertible and
\begingroup\postdisplaypenalty=100
\begin{align}
 &|\Sch_\fQ(\cM)(p,q)-\cM(p,q)|
 \le2\mathrm C_{3}a_*^2\Delta^{-1}M^{-r}\langle p-q\rangle^{-r},\quad p,q\in \fP,\ p\ne q,\label{2026102101}\\
 &|\Sch_\fQ(\cM)(p,p)-\cM(p,p)|
 \le2\mathrm C_{3}a_*^2\Delta^{-1}M^{-2r},\quad p\in \fP,\label{2026102102}\\
 &|\cW(\cM)(q,p)|
 \le2a_*\Delta^{-1}\langle q-p\rangle^{-r},\quad q\in \fQ,\ p\in \fP.\label{2026102103}
\end{align}
\endgroup
Here $\Sch_\fQ$ and $\cW(\cM)$ are defined in
\eqref{2026092803} and \eqref{2026092804}, respectively.
For $0\le\rho<r-d$, we also have
\begin{equation}\label{2026102104}
 \matnorm{\cW(\cM)}{\rho}
 \le1+\frac{\mathrm C_{4}(d,r,\rho)}{4}
       a_*\Delta^{-1}M^{\rho-r}.
\end{equation}
\end{lemma}
\begin{proof}
Let $\cD=\operatorname{diag}_{q\in \fQ}\cM(q,q)$ and
$\cM^{\mathrm{off}}=\cM_{\fQ}-\cD$.
By \eqref{2026101904} and \eqref{2026102004},
\[
 \sup_{q\in\fQ}\sum_{x\in\fQ}
       |(\cD^{-1}\cM^{\mathrm{off}})(q,x)|
 \le\mathrm C_{3} a_*\Delta^{-1}M^{-r}\le\frac12.
\]
It follows that
\[
 \cM_{\fQ}^{-1}
 =\sum_{\ell=0}^{\infty}(-\cD^{-1}\cM^{\mathrm{off}})^\ell\cD^{-1}.
\]
Multiplication by $\cD^{-1}\cM^{\mathrm{off}}$ involves a sum over
$x\in\fQ\setminus\{q\}$.
For $p\in\fP$ and $q\in\fQ$, \eqref{2026101905} gives
\begin{equation}\label{2026102105}
 \sum_{x\in\fQ\setminus\{q\}}
   |(\cD^{-1}\cM^{\mathrm{off}})(q,x)|\langle x-p\rangle^{-r}
 \le\mathrm C_{3} a_*\Delta^{-1}M^{-r}\langle q-p\rangle^{-r}
 \le\frac12\langle q-p\rangle^{-r}.
\end{equation}
Starting from
$|[\cD^{-1}\cM_{\fQ,p}](q)|\le a_*\Delta^{-1}\langle q-p\rangle^{-r}$,
induction using \eqref{2026102105} gives
\[
 |[(\cD^{-1}\cM^{\mathrm{off}})^\ell\cD^{-1}\cM_{\fQ,p}](q)|
 \le 2^{-\ell}a_*\Delta^{-1}\langle q-p\rangle^{-r},
 \qquad q\in\fQ,\quad\ell\ge0.
\]
After summing in $\ell$, we obtain
$|\cW(\cM)(q,p)|\le2a_*\Delta^{-1}\langle q-p\rangle^{-r}$, which proves
\eqref{2026102103}.
By \eqref{2026092803} and \eqref{2026102103},
\[
 |\Sch_\fQ(\cM)(p,p')-\cM(p,p')|
 \le2a_*^2\Delta^{-1}\sum_{q\in\fQ}
      \langle p-q\rangle^{-r}\langle q-p'\rangle^{-r},
 \qquad p,p'\in\fP.
\]
For $p\ne p'$, \eqref{2026101905} proves \eqref{2026102101}.
For $p=p'$, \eqref{2026101904} with $s=2r$ proves
\eqref{2026102102}.
Finally, multiply \eqref{2026102103} by
$\langle q-p\rangle^\rho$ and sum over $q$.
Since $\rho<r-d$, \eqref{2026101904} with $s=r-\rho$ proves
\eqref{2026102104} for a finite coefficient
$\mathrm C_4(d,r,\rho)$, as fixed before \eqref{2026100307}.
\end{proof}


\section{Geometric and domain verifications}\label{app:technical-verifications}

\subsection{The finite-region construction}\label{app:block-verification}
We prove Proposition~\ref{prop:det-blocks}, using the construction
\eqref{2026101503}--\eqref{2026101504}.
Fix $n\ge1$ and assume the geometric assertion \textup{(1)} of
Theorem~\ref{thm:finite-volume-induction} through step $n$.
\begin{proof}
Since \eqref{2026100703} gives
$|p-q|>33600N_h>18N_{h+1}$ for distinct $p,q\in Z_h$,
the balls $p+\Lambda_{9N_{h+1}}$, $p\in Z_h$, are pairwise disjoint.
For finite $B$, only finitely many centers $p\in Z_h$ satisfy
$(p+\Lambda_{9N_{h+1}})\cap B\ne\varnothing$.
Moreover, \eqref{2026101503} implies
\begin{equation*}
 B\subset\widetilde B\subset B+\Lambda_{18N_{h+1}},\qquad (p+\Lambda_{9N_{h+1}})\cap\widetilde B\ne\varnothing
 \quad\Longrightarrow\quad p+\Lambda_{9N_{h+1}}\subset\widetilde B,
 \qquad p\in Z_h.
\end{equation*}
Consequently,
\[
 \Lambda_{4N_{n+1}}\subset B_{n+1}
 \subset\Lambda_{4N_{n+1}+18\sum_{\ell=1}^{n}N_\ell}
 \subset\Lambda_{(26/5)N_{n+1}}\subset\Lambda_{8N_{n+1}}.
\]
Since $B_n\subset\Lambda_{8N_n}$ and $8N_n<4N_{n+1}$,
we obtain \eqref{2026101004} with $m=n+1$.

Let $0\le h\le n-1$ and $p\in B_{n+1}\cap Z_h$.
Let $F$ be the set immediately after the enlargement with centers in $Z_h$.
The subsequent enlargements give
\[
 B_{n+1}\subset F+\Lambda_{18\sum_{\ell=1}^{h}N_\ell},
 \qquad 18\sum_{\ell=1}^{h}N_\ell\le\frac65N_{h+1}<9N_{h+1}.
\]
Thus $p+\Lambda_{9N_{h+1}}$ meets $F$ and is contained in $F$.
By \eqref{2026101004}, $B_{h+1}\subset\Lambda_{8N_{h+1}}$, so
\[
 B_{h+1}(p)+\Lambda_{N_{h+1}}
 \subset p+\Lambda_{9N_{h+1}}\subset F\subset B_{n+1}.
\]
This proves \eqref{2026101005} for $m=n+1$.
If $p\in Z_n^\circ$, then
$|p|>33600N_n>9N_{n+1}$ by \eqref{2026100703}.
Since $B_{n+1}\subset\Lambda_{8N_{n+1}}$, we have
$B_{n+1}\cap Z_n=\{o\}$ and
$o\notin B_{n+1}(p)+\Lambda_{N_{n+1}}$.
Thus \eqref{2026101006} holds with $m=n+1$, and
\eqref{2026101007} holds with $h=n$.
\end{proof}

\subsection{Preservation of the multiplication domain}\label{app:domain-preservation}
Fix $j\ge1$ and assume assertions \textup{(1)--(4)} of
Theorem~\ref{thm:finite-volume-induction} through step $j+1$,
as in Proposition~\ref{prop:finite-vector-increment}.
Fix $\theta\in\Theta$, $t\in[0,1]$ and $i\in\{j,j+1\}$.
Use the matrices in \eqref{2026101607}--
\eqref{2026101608} and the multiplication
operator $\mathcal F$ in \eqref{2026101705}.
We prove \eqref{2026101706} and extend
\eqref{2026101701} from finitely supported
vectors to $\operatorname{Dom}(\mathcal F)$.
\begin{proof}
Since $\theta\in\Theta$, every diagonal entry of $\mathcal F$ is finite.
Each column $\mathcal V^{(i)}(\cdot,k)$ is supported on the finite set
$B_i(k)$, so $\mathcal V^{(i)}$ maps finitely supported vectors to
finitely supported vectors.
By \eqref{2026101102}, $\mathcal D^{(i)}-\mathcal F$ is bounded
on $\ell^1$.
Equation \eqref{2026101701} gives
\[
 [\mathcal F,\mathcal V^{(i)}]
 =\mathcal V^{(i)}(\mathcal D^{(i)}-\mathcal F)
    +\mathcal R^{(i)}-t\eps\cT\mathcal V^{(i)}
\]
on finitely supported vectors. The right-hand side is bounded on
$\ell^1$ by \eqref{2026101702},
\eqref{2026101704} and the boundedness of $\cT$.
For $u\in\operatorname{Dom}(\mathcal F)$, truncation to
$u_R=\Pi_{\Lambda_R}^*\Pi_{\Lambda_R}u$ yields
\begin{equation*}
 u_R\longrightarrow u,\qquad
       \mathcal F u_R\longrightarrow\mathcal F u,\qquad \mathcal F\mathcal V^{(i)}u_R
 \longrightarrow\mathcal V^{(i)}\mathcal F u
                    +[\mathcal F,\mathcal V^{(i)}]u
 \quad\text{in }\ell^1.
\end{equation*}
Since $\mathcal F$ is closed,
$\mathcal V^{(i)}u\in\operatorname{Dom}(\mathcal F)$.
Passing to the limit along $u_R$ therefore proves
\eqref{2026101701} for every
$u\in\operatorname{Dom}(\mathcal F)$.
For $\mathcal K^{(i)}=\mathcal V^{(i)}-\Id$ and $m\ge1$, the
commutator identity for $(\mathcal K^{(i)})^m$ gives
\begin{equation*}
 \|\mathcal F(\mathcal K^{(i)})^m u\|_1 \le(\mathrm C_{6} a)^m\|\mathcal F u\|_1
   +m(\mathrm C_{6} a)^{m-1}
     \|[\mathcal F,\mathcal V^{(i)}]\|_{\ell^1\to\ell^1}\|u\|_1.
\end{equation*}
Consequently, $\sum_{m\ge0}(-\mathcal K^{(i)})^m u$ converges
in the graph norm of $\mathcal F$.
Thus $(\mathcal V^{(i)})^{-1}$ also preserves
$\operatorname{Dom}(\mathcal F)$, which proves
\eqref{2026101706}.
\end{proof}

\raggedbottom

\end{document}